\documentstyle[amscd,12pt]{amsart}
\numberwithin{equation}{section}

\newcommand{\hb}{\hfill\break}
\newcommand{\eps}{\varepsilon}
\newcommand{\isum}{\mathop{\sum}_{i\in I}}
\newcommand{\hsum}{\mathop{\sum}_{\tau \in H}}

\newcommand{\iopl}{\mathop{\oplus}_{i\in I}}
\newcommand{\hopl}{\mathop{\oplus}_{\tau \in H}}
\newcommand{\omopl}{\mathop{\oplus}_{\tau \in \Omega}}
\def\boxtimes{\hbox{$\kern.5pt\sqcup\kern-8.0pt\sqcap$}
\kern-8.5pt\raise.6pt\hbox{$\times$}}
\newcommand{\om}{\Omega}
\newcommand{\evw}{E_{V,\Omega}}
\newcommand{\ev}{E_{V,\Omega}}
\newcommand{\xvw}{X_{V}}
\newcommand{\xv}{X_{V}}

\newcommand{\vin}{V_{\operatorname{in}(\tau)}}
\newcommand{\vout}{V_{\operatorname{out}(\tau)}}
\newcommand{\lvw}{\Lambda_{V}}

\newcommand{\lv}{\Lambda_{V}}

\newcommand{\D}{\cal D}
\newcommand{\pv}{{\cal P}_{V,\Omega}}
\newcommand{\qv}{{\cal Q}_{V,\Omega}}

\newcommand{\CO}{{\cal O}}
\newcommand{\cof}{{\operatorname{cof}\,}}
\renewcommand{\mod}{{\operatorname{mod}\,}}
\newcommand{\eit}{\tilde{e_i}}
\newcommand{\fit}{\tilde{f_i}}
\newcommand{\tf}{\tilde{f}}
\newcommand{\ft}{\tilde{f}}
\newcommand{\bb}{\bar{B}}

\newcommand{\dimc}{\dim}
\newcommand{\C}{{\Bbb C}}
\newcommand{\Z}{{\Bbb Z}}
\newcommand{\Homc}{\operatorname{Hom}\,}
\newcommand{\Hom}{\operatorname{Hom}\,}
\newcommand{\inn}{{\operatorname{in}}}
\newcommand{\tr}{{\operatorname{tr}}}
\newcommand{\out}{{\operatorname{out}}}
\newcommand{\Imm}{{\operatorname{Im}\,}}
\newcommand{\End}{{\operatorname{End}\,}}
\newcommand{\Ker}{{\operatorname{Ker}\,}}
\newcommand{\Coker}{{\operatorname{Coker}\,}}

\newcommand{\wt}{{\operatorname{wt}}}
\newcommand{\rank}{\operatorname{rank}}
\newcommand{\Gr}{\operatorname{Gr}}
\newcommand{\Ch}{\operatorname{Ch}}
\newcommand{\ran}{\rangle}
\newcommand{\lan}{\langle}
\newcommand{\gge}{>\kern-3pt>}

\newcommand{\te}{\tilde e}
\def\beq{\begin{eqnarray}}
\def\endeq{\end{eqnarray}}
\def\beqn{\begin{eqnarray*}}
\def\endeqn{\end{eqnarray*}}
\renewcommand\Im{\operatorname{Im}\,}
\theoremstyle{plain}
\newtheorem{thm}{Theorem}[subsection]
\newtheorem{prop}[thm]{Proposition}
\newtheorem{lemma}[thm]{Lemma}

\theoremstyle{definition}
\newtheorem{defn}[thm]{Definition}
\newtheorem{prob}{Problem}
\newtheorem{conj}[prob]{Conjecture}
\newtheorem{example}{Example}[subsection]
\newtheorem{rem}{Remark}

\renewcommand\max{{\operatorname{{max}}}\,}

\newenvironment{namelist}[1]{%
\begin{list}{}
  {
  \settowidth{\labelwidth}{#1}
  \setlength{\leftmargin}{1.1\labelwidth}}
}{%
\end{list}}
\begin{document}



\title[Geometric Construction of Crystal Bases]
{Geometric Construction of Crystal Bases}
\author{Masaki Kashiwara and Yoshihisa Saito}
\address[ ]{Research Institute for Mathematical Sciences,
Kyoto University, Japan}
\thanks{The second author is supported by the JSPS Research Fellowships for
Young Scientist.}
\maketitle
\begin{abstract}
We realize the crystal associated to
the quantized enveloping algebras
with a symmetric generalized Cartan matrix
as a set of
Lagrangian subvarieties of the cotangent bundle of the quiver variety.
As a by-product, 
we give a counterexample to the conjecture
of Kazhdan--Lusztig on the irreducibility
of the characteristic variety of the intersection cohomology sheaves 
associated with the Schubert cells of type A and also to
the similar problem asked by Lusztig on the characteristic variety of
the perverse sheaves corresponding to canonical bases.
\end{abstract}
\section{Introduction}
\subsection{}
G.Lusztig \cite{L;canI} gave a realization of
the quantized universal enveloping algebras
as the Grothendieck group of a category
of perverse sheaves on the quiver variety.
Let $(I,\Omega)$ be a finite oriented graph (=quiver), where
$I$ is the set of vertices and $\Omega$
is the set of arrows. Let us associate
a complex vector space $V_i$
to each vertex $i\in I$.
We set
$$ \evw =\omopl\Homc(V_{\out(\tau)},V_{\inn(\tau)})$$
and
$$ \xvw =\evw\oplus \evw^*
.$$
They are finite-dimensional vector spaces with
the action of the algebraic group $G_V=\prod_{i\in I}GL(V_i)$.
We regard $\xvw$ as the 
cotangent bundle of $\evw$.
Lusztig \cite{L;canI} realized a half of the quantized
universal enveloping algebra $U_q^-({\frak g})$ as the Grothendieck group of 
${\cal Q}_{V,\om}$. Here ${\cal Q}_{V,\om}$ is a subcategory of 
the derived category $D_c^b(E_{V,\om})$ 
of the bounded complex of constructible sheaves on 
$E_{V,\om}$. 
The irreducible perverse sheaves in  ${\cal Q}_{V,\om}$
form a base of $U_q^-({\frak g})$,
which is called canonical basis.

In \cite{L;quiver} he asked the following problem.\hb
{\bf Problem 1.}\quad If the underlying graph is of type $A$, $D$ or $E$,
then the singular support of any canonical base is irreducible.

One of the  purpose of this paper is to construct a counterexample of 
this problem for type $A$.
\subsection{} 
Let $G$ be a connected complex semisimple algebraic group,
$B$ a Borel subgroup of $G$ and $X=G/B$ the flag variety.
Let $D_X$ denote the sheaf of differential operators on $X$.
We denote the half sum of positive roots by $\rho$ 
and the Weyl group by $W$. For $w\in W$, let
$M_w$ be the Verma module with highest weight $-w(\rho)-\rho$ and $L_w$ its
simple quotient. By the Beilinson-Bernstein correspondence,
$M_w$ and $L_w$ correspond to regular holonomic $D_X$-modules
${\frak M}_w$ and ${\frak L}_w$ on $X$, respectively.
The characteristic varieties $\Ch({\frak M}_w)$ and 
$\text{Ch}({\frak L}_w)$ are Lagrangian subvarieties
of the cotangent bundle $T^*X$.
Each irreducible component of $\text{Ch}({\frak M}_w)$ and 
$\text{Ch}({\frak L}_w)$ is the closure of the conormal bundle $T^*_{X_y}X$ 
of a Schubert cell $X_y=ByB/B$ for some $y\in W$. Let ${\cal M}$ be the 
abelian category consisting of  regular holonomic systems on $X$ whose 
characteristic varieties are  contained in $\coprod_{w\in W}T^*_{X_w}X$. 
Its Grothendieck group  $K({\cal M})$ has two bases,
$([{\frak M}_w])_{w\in W}$ and $([{\frak L}_w])_{w\in W}$.
For ${\frak M}\in {\cal M}$ let $\text{\bf Ch}
({\frak M})=\sum_{w\in W}m_w({\frak M}) [T^*_{X_w}X]$ be the characteristic 
cycle. Here $m_w({\frak M})$ is the 
multiplicity of ${\frak M}$
along $T^*_{X_w}X$. Then
$\text{\bf Ch}$ extends to an additive map from $K({\cal M})$ to the group
of algebraic cycles of $T^*X$. Let $\chi$ be a $\Bbb Z$-linear isomorphism
from $K({\cal M})$ onto the group ring ${\Bbb Z}[W]$ defined by 
$\chi([{\frak M}_w])=w$. Then there exists a unique basis 
$\{\text{\bf b}(w)\}_{w\in W}$ of ${\Bbb Z}[W]$ such that $\text{\bf Ch}
(\chi^{-1}(\text{\bf b}(w)))=[T^*_{X_w}X]$ (See [KL1] and [KT].). This 
basis is related to the Springer representation of the Weyl group.
Set $\text{\bf a}(w)=\chi([{\frak L}_w])=\sum_{y\in W}m_y({\frak L}_w)
\text{\bf b}(y)$. The basis $\{\text{\bf a}(w)\}_{w\in W}$ is related to
the left cell representation of the Weyl group. 
Therefore an explicit knowledge of $m_y({\frak L}_w)$ gives
an explicit relation 
between the Springer representation and the left cell representation. 
If $\text{Ch}({\frak L}_w)$ is an irreducible variety, that is,
$$m_y({\frak L}_w)=
\begin{cases}
1 & \text{ if $y=w$,} \\
0 & \text{ otherwise},
\end{cases}
\leqno(1.2.1)
$$
then the Springer representation coincides with the left cell representation.
Due to Tanisaki, there is a counterexample of (1.2.1) in the case of 
$B_2$ (See [T]). In [KL2] Kazhdan and Lusztig conjectured that 
$\text{Ch}({\frak L}_w)$ is irreducible for $G=SL_n(\Bbb C)$. In this paper,
as a corollary of Problem 1, we shall show that there is a counterexample 
of this conjecture in the case of $G=SL_8(\Bbb C)$ and this conjecture
is true for $G=SL_n(\Bbb C)$ with
$n\leq 7$. 
\subsection{}
On the other hand, the first author [K1] constructed the crystal base and
the global crystal base of $U_q^-({\frak g})$
and the highest weight integrable representations of
$U_q({\frak g})$ in an algebraic way.
Grojnowski and Lusztig [GL]
showed that the global crystal base coincides with the
canonical base of Lusztig \cite{L;canI}. 

In this paper, we shall construct the crystal base in a 
geometrical way. We define the nilpotent subvariety
of the cotangent bundle of the quiver varieties, following Lusztig.
The nilpotent variety is a Lagrangian subvariety.
We shall define a crystal structure on the set
of its irreducible components, and we prove that it is
isomorphic to the crystal associated with $U_q^-({\frak g})$.

\subsection{}
 Let us briefly summarize the contents of this paper. In section 2 and 3 we 
 give a review of the theory of crystal base [K1,2,3,4].   
After recalling quiver varieties in section 4, we
define the crystal structure on the set of irreducible components
of the nilpotent varieties and prove that
it coincides with the crystal base of $U_q^-({\frak g})$
in section 5. In section 6, we recall the relation of 
the quantized universal enveloping
algebras and perverse sheaves on the quiver varieties.
In section 7, we give a negative answer to Problem 1.
In the last section, we give a counterexample of
the irreducibility of the characteristic variety of
the irreducible perverse sheaf with the Schubert cell as its support
in the case of $SL_8$.

We thank T.Tanisaki for stimulating discussions.

\section{Preliminaries}
\subsection{Definition of $U_q({\frak g})$}
We shall give the definition of $U_q({\frak g})$
associated with a symmetrizable Kac-Moody Lie algebra $\frak g$. 
We follow the notations in [K1,2,3,4].  

\begin{defn}
Let us consider following data:
\begin{namelist}{xxxx}
\item[(1)] a finite-dimensional $\Bbb Q$-vector space $\frak t$,
\item[(2)] an index set $I$ (of simple roots),
\item[(3)] a linearly independent subset $\{\alpha_i\,;\,i\in I\}$ of 
             ${\frak t}^*$ and a subset $\{h_i\,;\,i\in I\}$ of ${\frak t}$,
\item[(4)] an inner product ( , ) on ${\frak t}^*$ and
\item[(5)] a lattice $P$ (a weight lattice) of ${\frak t}^*$.
\end{namelist}
These data are assumed to satisfy the following conditions:
\begin{namelist}{xxxx}
\item[(6)] $\{\langle h_i,\alpha_j\rangle\}$ is a generalized Cartan 
          matrix\\
          (i.e. $\langle h_i,\alpha_i\rangle=2$, $\langle h_i,\alpha_j
          \rangle \in {\Bbb Z}_{\leq 0}$ for $i\ne j$ and 
          $\langle h_i,\alpha_j\rangle =0\Leftrightarrow 
          \langle h_j,\alpha_i\rangle =0$),
\item[(7)] $(\alpha_i,\alpha_i)\in {2\Bbb Z}_{>0},$
\item[(8)] $\langle h_i,\lambda\rangle =2(\alpha_i,\lambda)/
         (\alpha_i,\alpha_i)$
             for any $i\in I$ and $\lambda \in {\frak t}^*$,
\item[(9)] $\alpha_i\in P$ and $h_i\in P^*=\{h\in {\frak t}~\,;\,
         \langle h,P\rangle\in {\Bbb Z}\}$.
\end{namelist}
Then the ${\Bbb Q}(q)$-algebra $U_q({\frak g})$ is the algebra generated by
$e_i,f_i(i\in I)$ and $q^h(h\in P^*)$ with the following defining relations:
\begin{namelist}{xxxx}
\item[(10)] $q^h=1$ for $h=0$ and $q^{h+h'}=q^hq^{h'}$,
\item[(11)] $q^he_iq^{-h}=q^{\langle h,\alpha_i\rangle}e_i$ and 
          $q^hf_iq^{-h}=q^{-\langle h,\alpha_i\rangle}f_i$,
\item[(12)] $[e_i,f_j]=\delta_{i,j}(t_i-t_i^{-1})/(q_i-q_i^{-1})$
              where $q_i=q^{(\alpha_i,\alpha_i)/2}$ and $t_i=q^{(\alpha_i,
              \alpha_i)h_i/2}$,
\item[(13)] $\sum\limits_{n=0}^{b}(-1)^ne_i^{(n)}e_je_i^{(b-n)}=
\sum\limits_{n=0}^{b}(-1)^nf_i^{(n)}f_j
              f_i^{(b-n)}=0$\\
              where $i\ne j$ and $b=1-\langle h_i,\alpha_j\rangle$.
\end{namelist}
\end{defn}
Here we used the notations $[n]_i=(q_i^n-q_i^{-n})/(q_i-q_i^{-1})$, $[n]_i!
=\prod_{k=1}^n [k]_i$, $e_i^{(n)}=e_i^n/[n]_i!$ and $f_i^{(n)}=f_i^n/[n]_i!$.
We understand $e_i^{(n)}=f_i^{(n)}=0$ for $n<0$. We set
$Q=\isum {\Bbb Z}\alpha_i$, $Q_+=\isum {\Bbb Z}_{\geq 0}\alpha_i$
and $Q_-=-Q_+$. Let $P_+$ be the set of dominant integral weights.
\hb
We denote by $U_q^-({\frak g})$ the ${\Bbb Q}(q)$-subalgebra of 
$U_q({\frak g})$ generated by $f_i$ $(i\in I)$.

As in [K], we define the ${\Bbb Q}(q)$-algebra anti-automorphism $*$ of
$U_q({\frak g})$ by
$$\text{${e_i}^*=e_i$, ${f_i}^* =f_i$ and $(q^h)^*=q^{-h}$.}$$
Note that $*^2=1$.
\subsection{Crystal base}
In this subsection we give a review of the theory of crystal base. See 
[K1,2,3,4] for details.

Let $M$ be an integrable $U_q({\frak g})$-module and let $M=\oplus_{\nu\in P}
M_{\nu}$ be the weight space decomposition. By the theory of 
integrable representation of $U_q({\frak {sl}}(2))$, we have
$$M=\mathop{\oplus}_{0\leq n\leq \langle h_i,\nu\rangle} f_i^{(n)}
(\Ker e_i\cap M_{\nu}).$$
We define the endomorphisms $\eit$ and $\fit$ of $M$ by
$$\fit(f_i^{(n)}u)=f_i^{(n+1)}u \text{ and}$$
$$\eit(f_i^{(n)}u)=f_i^{(n-1)}u$$
for $u\in \Ker e_i\cap M_{\nu}$ with $0\leq n\leq \langle h_i,\nu\rangle$.

Let $A$ be the subring of ${\Bbb Q}(q)$ consisting of rational functions 
without pole at $q=0$.

\begin{defn}
A pair $(L,B)$ is called a crystal base of $M$ if it satisfies the following 
conditions:
\begin{namelist}{(2.2.1)xxxx}
\item[(2.2.1)] $L$ is a free sub-$A$-module of $M$ such that $M\cong
             {\Bbb Q}(q)\otimes_A L$.
\item[(2.2.2)] $B$ is a base of the $\Bbb Q$-vector space $L/qL$.
\item[(2.2.3)] $\eit L\subset L$ and $\fit L\subset L$ for any $i$.
\item[{\phantom{(2.2.3)}}] Therefore, $\eit$ and $\fit$ act on $L/qL$.
\item[(2.2.4)] $\eit B\subset B\sqcup \{0\}$ and $\fit B\subset B\sqcup 
             \{0\}$.
\item[(2.2.5)] $$L=\mathop{\oplus}_{\nu\in P}L_{\nu} \text{ and }
               B=\mathop{\bigsqcup}_{\nu\in P}B_{\nu}$$
             where $L_{\nu}=L\cap M_{\nu}$ and $B_{\nu}=B\cap 
             (L_{\nu}/qL_{\nu})$.
\item[(2.2.6)] For $b$, $b'\in B$, $b'=\fit b$ if and only if $b=\eit b'$.
\end{namelist}
\end{defn}

For \(\lambda\in P_+\), we denote by $V(\lambda)$ the simple
$U_q({\frak g})$-module of highest weight $\lambda$. The highest weight 
vector of $V(\lambda)$ is denoted by $u_{\lambda}$. We consider
the sub-$A$-module $L(\lambda)$ of $V(\lambda)$ generated by $\tilde{f_{i_1}}
\cdots\tilde{f_{i_l}}u_{\lambda}$ and the subset $B(\lambda)$ of
$L(\lambda)/qL(\lambda)$ consisting of the non-zero vectors of the form
$\tilde{f_{i_1}}\cdots\tilde{f_{i_l}}u_{\lambda}$.

\begin{thm}
$(L(\lambda),B(\lambda))$ is a crystal base of $V(\lambda)$. 
\end{thm}

\subsection{Crystal base of $U_q^-({\frak g})$}
Next we shall define a crystal base of $U_q^-({\frak g})$.

\begin{lemma}
For any $P\in U_q^-({\frak g})$, there exist unique $Q$, $R\in 
U_q^-({\frak g})$ such that
$$[e_i,P]=\frac{t_iQ-t_i^{-1}R}{q_i-q_i^{-1}}.$$
\end{lemma}
By this lemma, $e_i'(P)=R$ defines an endomorphism $e_i'$ of
$U_q^-({\frak g})$.

According to [K1] we have
$$U_q^-({\frak g})=\mathop{\oplus}_{n\geq 0} f_i^{(n)}\Ker e'_i.$$
We define the endomorphisms $\eit$ and $\fit$ of $U_q^-({\frak g})$ by
\beqn
\fit(f_i^{(n)}u)&=&f_i^{(n+1)}u\quad\text{and}\\
\eit(f_i^{(n)}u)&=&f_i^{(n-1)}u
\endeqn
for $u\in\Ker e'_i$.

\begin{defn}
A pair $(L,B)$ is called a crystal base of $U_q^-({\frak g})$ if it 
satisfies the following conditions:
\begin{namelist}{(2.3.1)xxxx}
\item[(2.3.1)] $L$ is a free sub-$A$-module of $U_q^-({\frak g})$ such that 
             $U_q^-({\frak g})\cong {\Bbb Q}(q)\otimes_A L$.
\item[(2.3.2)] $B$ is a base of the $\Bbb Q$-vector space $L/qL$.
\item[(2.3.3)] $\eit L\subset L$ and $\fit L\subset L$ for any $i$.
\item[{\phantom{(2.3.3)}}] Therefore $\eit$ and $\fit$ act on $L/qL$.
\item[(2.3.4)] $\eit B\subset B\sqcup \{0\}$ and $\fit B\subset B$.
\item[(2.3.5)] $$L=\mathop{\oplus}_{\nu\in Q_-}L_{\nu} \text{ and }
               B=\mathop{\bigsqcup}_{\nu\in Q_-}B_{\nu}$$
             where $L_{\nu}=L\cap U_q^-({\frak g})_{\nu}$, $B_{\nu}=B\cap 
             (L_{\nu}/qL_{\nu})$ and
$U_q^-({\frak g})_\nu
=\{P\in U_q^-({\frak g})\,;\,q^hPq^{-h}=q^{\lan h,\nu\ran}P
\text{ for any $h\in P^*$}\}$.
\item[(2.3.6)] For $b\in B$ such that $\eit b\not=0$,
we have $b=\fit\eit b$.
\end{namelist}
\end{defn}

We introduce the sub-$A$-module $L(\infty)$ of $U_q^-({\frak g})$ generated 
by $\tilde{f_{i_1}}\cdots\tilde{f_{i_l}}\cdot 1$ and the subset $B(\infty)$ 
of $L(\infty)/qL(\infty)$ consisting of the non-zero vectors of the form
$\tilde{f_{i_1}}\cdots\tilde{f_{i_l}}\cdot 1$.

\begin{thm}
$(L(\infty),B(\infty))$ is a crystal base of $U_q^-({\frak g})$.
\end{thm}

\section{Crystals} 
\subsection{Definition of Crystal} 
\begin{defn}\label{dfn:cry}
A crystal B is a set endowed with
$$\text{maps }\ \wt:B\to P,\mbox{ } \varepsilon_i:B\to {\bf Z}\sqcup
\{-\infty \},\mbox{ } \varphi_i:B\to {\bf Z}\sqcup \{-\infty \}\ 
\mbox{ and}\leqno{(3.1.1)}$$
$$\tilde{e}_i :B\to B\sqcup \{0\},\mbox{ } \tilde{f}_i :B\to B\sqcup \{0\}. 
\leqno{(3.1.2)}$$
They are subject to the following axioms:
\begin{namelist}{(C 1)xxxx}
\item[(C 1)] \(\varphi_i (b)=\varepsilon_i (b)+\langle h_i,\wt(b)\rangle.\)
\item[(C 2)] If \(b\in B\) and \(\tilde{e}_i b\in B\) then,\hb
\(\wt(\tilde{e}_i b)=\wt(b)+\alpha_i\), \(\varepsilon_i(\tilde{e_i} b)=
\varepsilon_i (b)-1\) and \(\varphi_i (\tilde{e_i} b)=\varphi_i(b)+1.\)
\item[(C 2')] If \(b\in B\) and \(\tilde{f_i} b\in B\), then\hb
\(\wt(\tilde{f_i} b)=\wt(b)-\alpha_i\), \(\varepsilon_i (\tilde{f_i} b)=
\varepsilon_i(b)+1\) and \(\varphi_i(\tilde{f_i} b)=\varphi_i(b)-1\).
\item[(C 3)] For \(b,b'\in B\) and \(i\in I\), \(b'=\tilde{e_i} b\) if and 
only if
\(b=\tilde{f_i} b'\).
\item[(C 4)] For \(b\in B\), if \(\varphi_i(b)=-\infty \), then 
\(\tilde{e_i} b
=\tilde{f_i} b=0.\)
\end{namelist}
\end{defn}
For two crystals \(B_1\) and \(B_2\), a morphism \(\psi\) from \(B_1\) to
\(B_2\) is a map \(B_1\sqcup \{0\}
\to B_2\sqcup \{0\}\) that satisfies the following
conditions:
\begin{namelist}{(3.1.3)xxxx}
\item[(3.1.3)] $\psi(0)=0$,
\item[(3.1.4)] If \(b\in B_1\) and \(\psi (b)\in B_2\), then\hb
             \(\wt(\psi(b))=\wt(b)\), \(\varepsilon_i (\psi (b))=
             \varepsilon_i (b)\), and \(\varphi_i(\psi (b))=\varphi_i 
             (b)\),
\item[(3.1.5)] If \(b,b'\in B_1\) and $i\in I$ satisfy
$\fit(b)=b'$ and $\psi(b)$, $\psi(b')\in B_2$, 
then we have $\fit(\psi(b))=\psi(b')$.
\end{namelist}
A morphism \(\psi :B_1\to B_2\) is called {\it strict}, if it commutes with
all \(\tilde{e_i}\) and \(\tilde{f_i}\).

A morphism \(\psi :B_1\to B_2\) is called an {\it embedding}, if \(\psi\)
induces an injective map from \(B_1\sqcup \{0\}\) to \(B_2\sqcup \{0\}\).

For two crystals \(B_1\) and \(B_2\), we define its tensor product
\(B_1\otimes B_2\) as follows:
\beqn
B_1\otimes B_2&=&\left\{ 
b_1\otimes b_2~;~b_1\in B_1\mbox{ } and\mbox{ } b_2
\in B_2\right\},\\
\varepsilon_i (b_1\otimes b_2)&=&\max\Big(\varepsilon_i (b_1),\mbox{ }
\varepsilon_i (b_2)-\wt_i(b_1)\Big),\\
\varphi_i (b_1\otimes b_2)&=&\max\Big(\varphi_i (b_1)+\wt_i(b_2),\mbox{ }
\varphi_i (b_2)\Big),\\
\wt(b_1\otimes b_2)&=&\wt(b_1)+\wt(b_2).
\endeqn
Here \(\wt_i(b)\) denotes \(\langle h_i,\wt(b)\rangle\).

The action of \(\tilde{e_i}\) and \(\tilde{f_i}\)
are defined by
\[\tilde{e_i} (b_1\otimes b_2) = \left\{
\begin{array}{rl}
\tilde{e_i} b_1\otimes b_2 &\quad\mbox{if $\varphi_i (b_1)\geq \varepsilon_i
(b_2)$}\\
b_1\otimes \tilde{e_i} b_2 &\quad\mbox{if $\varphi_i (b_1)<\varepsilon_i
(b_2)$,}
\end{array}\right. \]
\[\tilde{f_i} (b_1\otimes b_2) = \left\{
\begin{array}{rl}
\tilde{f_i} b_1\otimes b_2 &\quad\mbox{if $\varphi_i (b_1)>\varepsilon_i
(b_2)$}\\
b_1\otimes \tilde{f_i} b_2 &\quad\mbox{if $\varphi_i (b_1)\leq \varepsilon_i
(b_2)$}.
\end{array}\right. \]

\begin{example}
For \(i\in I\), \(B_i\) is the crystal defined as follows
\beqn
&&B_i=\left\{b_i(n)\,;\,n\in {\Bbb Z}\right\},\\
&&\wt(b_i(n))=n\alpha_i,\\
&&\varphi_i (b_i(n))=n,\mbox{ }\varepsilon_i (b_i(n))=-n,\\
&&\varphi_j (b_i(n))=\varepsilon_j (b_i(n))=-\infty \quad
\text{for $i\ne j$}.
\endeqn
We define the action of $\eit$ and $\fit$ by
\beqn
&&\tilde{e_i} (b_i(n))=b_i(n+1),\\
&&\tilde{f_i} (b_i(n))=b_i(n-1),\\
&&\tilde{e_j} (b_i(n))=\tilde{f_j} (b_i(n))=0 \quad\mbox{for $i\ne j$}.
\endeqn
We write \(b_i\) for \(b_i(0)\).
\end{example}

\begin{example}
For \(\lambda\in P_+\), \(B(\lambda )\) denotes the crystal associated with the
crystal base of the simple highest weight module with highest weight
\(\lambda \). For $b\in B(\lambda)$ we set $\varepsilon_i (b)
=\max\{k\geq 0~;~\tilde{e_i}^k b\ne 0\}$, $\varphi_i (b)
=\max\{k\geq 0~;~\tilde{f_i}^k b\ne 0\}$ and $\wt(b)$ is the weight of $b$.
\end{example}

\begin{example}
\(B(\infty )\) is the crystal associated with the crystal base of 
$U_q^-({\frak g})$. For $b\in B(\infty)$ we set $\varepsilon_i (b)=
\max\{k\geq 0~;~\tilde{e_i}^k b\ne 0\}$ and $\varphi_i (b)=\varepsilon_i (b)
+\langle h_i,\wt(b)\rangle$. We denote $u_{\infty}$ by the unique element
with weight $0$.
\end{example}
\subsection{}
We have $L(\infty)^*=L(\infty)$ and
$*$ induces an endomorphism of $L(\infty)/qL(\infty)$.
\begin{thm}
$$B(\infty )^*=B(\infty ).$$
\end{thm}
We define the operators \(\tilde{e}_i^*\), \(\tilde{f}_i^*\) of
\(U_q^-({\frak g})\)
by
$$\tilde{e}_i^* =*\tilde{e}_i *,\mbox{ and }\tilde{f}_i^* =*\tilde f_i *.
\leqno{(3.2.2)} $$

\begin{thm}
\begin{enumerate}
\item For any i, there exists a unique strict embedding of crystals
\[\Psi_i :B(\infty )\hookrightarrow B(\infty )\otimes B_i\]
that sends \(u_{\infty }\) to \(u_{\infty }\otimes b_i\).
\item If \(\Psi_i (b)=b'\otimes \tilde{f}_i^n  b_i\mbox{ }(n\geq 0)\),
then $\varepsilon_i (b^*)=n$, $\varepsilon_i ({b'}^*)=0$ and
$b=\fit^*{}^nb'$.
\item $\Imm~\Psi_i =\{b\otimes \tilde{f_i}^{n}b_i~\,;\,b\in B(\infty ),
\mbox{ }\varepsilon_i (b^*)=0,\mbox{ }n\geq 0\}$.
\end{enumerate}
\end{thm}

%

In fact the above properties characterize
$B(\infty)$ as seen in the following proposition.

\begin{prop}\label{char}
Let $B$ be a crystal and $b_0$ an element of $B$ with weight $0$.
Assume the following conditions.
\begin{enumerate}
\item\label{pr:1}
$\wt(B)\subset Q_-$.
\item\label{pr:2}
$b_0$ is a unique element of $B$ with weight $0$.
\item
$\eps_i(b_0)=0$ for every $i$.
\item
$\eps_i(b)\in\Z$ for any $b$ and $i$.
\item
For every $i$, there exists
a strict embedding $\Psi_i:B\to B\otimes B_i$.
\item\label{pr:3}
$\Psi_i(B)\subset B\times\{\fit^nb_i;n\ge0\}$.
\item\label{cond:1}
For any $b\in B$ such that $b\not=b_0$, there exists $i$ such that
$\Psi_i(b)=b'\otimes \fit^nb_i$ with $n>0$.
\end{enumerate}
Then $B$ is isomorphic to $B(\infty)$.
\end{prop}

\begin{pf}
First note that
$\Psi_i(b_0)=b_0\otimes b_i$
by (\ref{pr:1}), (\ref{pr:2}) and (\ref{pr:3}).

We shall show 
that for any $b\in B$ with $b\not=b_0$
there exists $i$ such that $\eit(b)\not=0$.
Take $i$ such that
$\Psi_i(b)=b'\otimes \fit^nb_i$ with $n>0$.
If $b'=b_0$ then
$\eit(b)=b_0\otimes \fit^{n-1}b_i\not=0$.
If $b'\not=b_0$, then
the induction on the weight implies
the existence of $j\in I$ such that $\tilde{e_j}(b')\not=0$.
Then $\tilde{e_j}(b)\not=0$.

Hence any element of $B$ has the form
$\tf_{i_1}\cdots\tf_{i_l}b_0$ with $i_1,\cdots i_l\in I$.

Now take a sequence $(i_1,i_2,\cdots)$ in $I$ in which
every element of $I$ appears infinitely many times.
Let us consider the composition $\Phi_n$ of the following
chain of crystal morphisms.
\beqn
\Phi_n:B&@>{\Psi_{i_1}}>>&B\otimes B_{i_1}
@>{\Psi_{i_2}}>>B\otimes B_{i_2}\otimes B_{i_1}@>{}>>\cdots
@>{\Psi_{i_n}}>>B\otimes B_{i_n}\otimes\cdots\otimes B_{i_1}.
\endeqn
Then for any $b\in B$ there exists $n$ such that
$\Phi_n(b)$ has the form
$b_0\otimes\tilde{f}_{i_n}^{a_n}b_{i_n}
\otimes\tilde{f}_{i_{n-1}}^{a_{n-1}}b_{i_{n-1}}
\otimes\cdots\otimes\tilde{f}_{i_1}^{a_1}b_{i_1}$.
The sequence $(a_1,a_2,\cdots,a_n,0,0,\cdots)$
does not depend on such a choice of $n$.
Let $\tilde B$ be the set of sequences
$(a_1,a_2,\cdots,a_n,0,0,\cdots)$ of 
integers such that $a_n=0$ for $n\gge0$.
Then $\tilde B$ has a crystal structure
by
$(a_1,a_2,\cdots,a_n,0,0,\cdots)\mapsto \cdots\otimes
b_{i_{2}}(-a_2)\otimes b_{i_1}(-a_1)$.
Then both $B(\infty)$ and $B$ are strictly embedded into
$\tilde B$ and their images coincide
with the smallest strict subcrystal of $\tilde B$
containing $(0,0,\cdots)$.
Therefore, they are isomorphic.
\end{pf}


\section{Quivers and associated varieties ([L5,6] and [N1,2])}
\subsection{Definition of quiver}
We shall recall the formulation due to Lusztig [L5,6].

Suppose a finite graph is given. In this graph, two different 
vertices may be joined by several edges, 
but any vertex is not joined with itself
by any edges.
Let $I$ be the set of vertices of our graph, and let $H$
be the  set of pairs of an edge and its orientation.
The precise definition is as follows.

\begin{defn}
Suppose that following data (1) $\sim$ (5) are given:
\begin{enumerate}
\item  a finite set $I$,  
\item a finite set $H$,
\item a map $H \to I$ denoted $\tau \mapsto \out(\tau)$,
\item a map $H \to I$ denoted $\tau \mapsto \inn(\tau)$ and 
\item an involution $\tau \mapsto \bar{\tau}$ of $H$.
\end{enumerate}
We assume that they satisfy the following conditions;
$$ \inn(\bar{\tau})=\out(\tau),\ \out(\bar{\tau})=\inn(\tau)\ \text{ and} 
\leqno{(4.1.1)} $$
$$ \out(\tau)\ne \inn(\tau) \text{ for all $\tau\in H$}.\leqno{(4.1.2)} $$
An orientation of the graph is a choice of a subset $\Omega \subset H$ such 
that
$$ \Omega \cup \bar{\Omega} =H \text{ and } 
\Omega \cap \bar{\Omega} =\phi.$$
We call a quiver a graph with an orientation.
\end{defn}
To a graph $(I,H)$ we associate a root system
with simple roots $\{\alpha_i\}_{i\in I}$
and simple coroots $\{h_i\}_{i\in I}$
with
$$\langle h_i,\alpha_j\rangle=(\alpha_i,\alpha_j)=
\begin{cases}
2, & i=j,\\
-\sharp\{\tau\in H~;~\out(\tau)=i,~\inn(\tau)=j\}, & i\ne j.
\end{cases}
$$
We denote by $\frak g$ the corresponding
Kac-Moody Lie algebra and $U_q(\frak g)$ the corresponding
quantized universal enveloping algebra.

\subsection{ }
Let $\cal V$ be the family of $I$-graded complex vector spaces 
$V=\iopl V_i$.
We set $\dim V=-\sum_{i\in I}(\dim V_i)\alpha_i\in Q_-$.
For 
$\nu\in Q_-$,
let ${\cal V}_{\nu}$ be the family of 
$I$-graded complex vector spaces $V$ 
with $\dim V=\nu$.

Let us define the complex vector spaces $\ev$ and $\xv$ by
\begin{eqnarray*}
\ev&=&\omopl \Homc( V_{\out(\tau)},
 V_{\inn(\tau)}),\\
 \xv& =&\hopl \Homc(\vout ,\vin).
\end{eqnarray*}

In the sequel, a point of $E_{V,\Omega}$ or $X_{V}$ will be denoted
as $B=(B_\tau)$. Here $B_\tau$ is in $\Homc ( V_{\out(\tau)},V_{\inn(\tau)})$.

We define the symplectic form $\omega $ on $X_{V}$ by
$$ \omega (B,B')=
 \hsum \eps (\tau)\tr(B_{\bar{\tau}}B'_{\tau}),\leqno(4.2.1)$$
where $\eps (\tau)=1$ if $\tau\in \om$, $\eps (\tau)=-1$ if $\tau\in 
\bar{\om}$. We sometimes identify $\xvw$ and the cotangent bundle of $\evw$ 
via $\omega$.

The group $G_V=\prod_{i\in I}GL(V_i)$ acts on $\evw$ and $\xvw$ by
\begin{eqnarray*}
G_V\ni g=(g_i)\,&:&\,
(B_{\tau})\mapsto(g_{\inn(\tau)}B_{\tau}g_{\out(\tau)}^{-1}),
\end{eqnarray*}
where $g_i\in GL(V_i)$ for each $i\in I$.

The Lie algebra of $G_V$ is ${\frak g}_V=\bigoplus_{i\in I}\End(V_i)$. 
We denote an element of ${\frak g}_V$ by 
$A=(A_i)_{i\in I}$ with $A_i\in\End(V_i)$. 
The infinitesimal action of $A\in {\frak g}_V$ on $\xvw$ 
at $B\in \xvw$ is given by
$[A,B]$.
Let  $\mu :\xvw \to {\frak g}_V$ be the moment map
associated with the $G_V$-action on the symplectic vector space
$\xvw$. Its $i$-th component $\mu_i:\xvw \to \End(V_i)$
is given by
\[\mu_i (B)=\sideset{}{}\sum
\begin{Sb}
\tau\in H \\ i=\out(\tau) 
\end{Sb} \eps(\tau)B_{\bar{\tau}}B_{\tau}.\]

For a non-negative integer $n$, we set 
$${\frak S}_n=\{ \sigma =(\tau_1,\tau_2,
\cdots,\tau_n)\,;\,\tau_i\in H,
\inn(\tau_1)=\out(\tau_2),\cdots,\inn(\tau_{n-1})=\out(\tau_n)\}\,,$$
 and set
${\frak S}=\bigcup_{n\ge 0}{\frak S}_n$.
For $\sigma =(\tau_1,\tau_2,\cdots,\tau_n)$, we
set $\out(\sigma)=\out(\tau_1),\inn(\sigma)=\inn(\tau_n)$.
For $B\in X_V$ we set
$B_{\sigma}=B_{\tau_n}\cdots B_{\tau_1}
:V_{\out(\tau_1)}\to V_{\inn(\tau_n)}$.
If $n=0$, we understand that
${\frak S}_n=\{1\}$ and $B_1$ is the identity.
An element $B$ of $X_V$ is called nilpotent if there exists a
positive integer $n$ such that $B_{\sigma}=0$
for any $\sigma\in{\frak S}_n$.

\begin{defn}
We set
$$X_0{}_V=\{B\in \xvw\,;\,\mu(B)=0\}$$
and
$$\lvw =\{ B\in \xvw\,;\,\mu(B)=0\text{ and } B\text{
is nilpotent}\}.$$
\end{defn}
It is clear that $\lvw$ is a $G_V$-stable 
closed subvariety of $\xvw$.
It is known that $\lvw$ is a Lagrangian variety \cite{L;quiver}.
\section{Lagrangian construction of crystal base}
\subsection{}
For each $\nu\in Q_-$, let us take $V(\nu)\in{\cal V}_\nu$
and set $X(\nu)=X_{V(\nu)}$,  $X_0(\nu)=X_0{}_{\,V(\nu)}$ and
$\Lambda(\nu)=\Lambda_{V(\nu)}$.
For $\nu$, $\nu'$ and $\bar{\nu}$ in $Q_-$ with $\nu=\nu'+\bar{\nu}$,
we consider the diagram
\begin{eqnarray}\label{seqx}
&&X_0({\bar\nu})\times X_0({\nu'})\mathop{\longleftarrow}^{q_1}X'_0(\bar\nu,\nu')
\mathop{\longrightarrow}^{q_2}
X_0(\nu).
\end{eqnarray}
Here
$X'_0(\bar\nu,\nu')$ is the variety of $(B,\bar{\phi},\phi')$ 
where $B\in X_0(\nu)$ and
$\bar{\phi}=(\bar\phi_i)$, $\phi=(\phi'_i)$ give an exact
sequence
\beq
0\to V(\bar\nu)_i\mathop{\to}^{\bar\phi_i}V(\nu)_i\mathop{\to}^{\phi'_i}
{V(\nu')}_i\to 0\label{dia:2}
\endeq
such that $\Imm\bar\phi$ is 
stable by $B$. Hence $B$ induces
$\bar{B}:{V(\bar\nu)}\to\bar{V(\bar\nu)}$ and $B':V(\nu')\to V(\nu')$.
We define  $q_1(B,\bar\phi,\phi')=(\bb,B')$ and
$q_2(B,\bar\phi,\phi')=B$.

The following lemma is easily proved.

\begin{lemma}
Under the above notations the following two conditions are equivalent.
     \begin{enumerate}
      \renewcommand{\labelenumi}{(\alph{enumi})}%
       \item $B$ is nilpotent.
       \item Both $B'$ and $\bb$ are nilpotent.
     \end{enumerate}
\end{lemma}

By this lemma, the diagram (\ref{seqx}) induces the diagram
\begin{eqnarray}\label{seql}
&&\Lambda(\bar\nu)\times \Lambda(\nu')\mathop{\longleftarrow}^{q_1}\Lambda'
(\bar\nu,\nu')\mathop{\longrightarrow}^{q_2}
\Lambda(\nu).
\end{eqnarray}
Here $\Lambda'(\bar\nu,\nu')=q_2^{-1}\big(\Lambda(\nu)\big)
=q_1^{-1}\big(\Lambda(\bar\nu)\times \Lambda(\nu')\big)$.

For $i\in I$ and $p\in {\Bbb Z}_{\geq 0}$ we consider
$$X_0(\nu)_{i,p}=\{B\in X_0(\nu)\,;\,\eps_i(B)=p\},$$
where
$$\eps_i(B)=\dimc\Coker\Big(\mathop{\oplus}_{\tau;\inn(\tau)=i}
V(\nu)_{\out(\tau)}@>{(B_{\tau})}>>V(\nu)_i\Big).$$
It is clear that $X_0(\nu)_{i,p}$ is a locally closed subvariety of 
$X_0(\nu)$.
\subsection{}
In this and the next subsections,
we assume that $\nu=\bar{\nu}-c\alpha_i$ 
for $c\in {\Bbb Z}_{\geq 0}$.
We set $V=V(\nu)$ and $\bar V=V(\bar\nu)$.

Let us consider the special case of (\ref{seqx}).
Note that $X_0(-c\alpha_i)=\{0\}$.
\begin{eqnarray}
&&X_0({\bar\nu})\cong X_0({\bar\nu})\times X_0({-c\alpha_i})
\mathop{\longleftarrow}^{\varpi_1}
X'_0(\bar\nu,-c\alpha_i)
\mathop{\longrightarrow}^{\varpi_2}
X_0(\nu).\label{dia:e}
\end{eqnarray}

\begin{lemma}
Let $p\in {\Bbb Z}_{\geq 0}$. Then we have
$$\varpi_1^{-1}(X_0({\bar\nu})_{i,p})=\varpi_2^{-1}
(X_0(\nu)_{i,\,p+c}).$$
\end{lemma}

\begin{pf}
This follows immediately from the diagram (\ref{dia:2}).
\end{pf}

\begin{defn}
We set
\begin{align*}
X'_0(\bar{\nu},-c\alpha_i)_p
& =\varpi_1^{-1}(X_0(\bar{\nu})_{i,p})
=\varpi_2^{-1}(X_0(\nu)_{i,\,p+c}).
\end{align*}  
\end{defn}

Suppose $p=0$. Then we have following diagram
\begin{eqnarray}
X_0(\bar{\nu})\supset X_0(\bar{\nu})_{i,0} 
\mathop{\longleftarrow}^{\varpi_1}
X'_0(\bar{\nu},-c\alpha_i)_0 \mathop{\longrightarrow}^{\varpi_2} X_0(\nu)_{i,c} 
\subset X_0(\nu).
\end{eqnarray}
Note that $X_0(\bar{\nu})_{i,0}$ is 
an open subvariety of $X_0(\bar{\nu})$.

\begin{lemma}
\begin{enumerate}
\item
$\varpi_2:X'_0(\bar{\nu},-c\alpha_i)_0 \mathop{\longrightarrow}
X_0(\nu)_{i,c} $
is a principal fiber bundle
with $GL(\C^c)\times
\prod_{j\in I}GL(V(\bar\nu)_j)$ as fiber.
\item
${\varpi_1}:X'_0(\bar{\nu},-c\alpha_i)_0
\mathop{\longrightarrow}X_0(\bar{\nu})_{i,0}$
is a smooth map whose fiber is a connected
rational variety of dimension
$\sum_{j\in I}\Big(\dim V(\nu)_j\Big)^2-c\,(\alpha_i,\bar\nu)$.
\end{enumerate}
\end{lemma}

\begin{pf}

\noindent
(1)\quad The fiber of $\varpi_2$ over $B\in X_0(\nu)_{i,c}$
is the set of families of isomorphisms 
\beqn
&&V(\bar\nu)_j{\buildrel \sim\over\to}V(\nu)_j\quad(j\not=i),\\
&&V(\bar\nu)_i{\buildrel \sim\over\to}
\Imm(\mathop{\oplus}_{\inn(\tau)=i}
{V}(\nu)_{\out(\tau)}\to V(\nu)_i)\quad\mbox{and}\\
&&\Coker(\mathop{\oplus}_{\inn(\tau)=i}
{V}(\nu)_{\out(\tau)}\to V(\nu)_i){\buildrel \sim\over\to}
V(-c\alpha_i)_i.
\endeqn
Hence we obtain (1).

\noindent
(2)\quad Let $\bb\in X_0(\bar{\nu})_{i,0}$ and take
$(B,\bar\phi,\phi')\in\varpi_1^{-1}(\bar B)$.
Consider the following diagram\,:
$$
\begin{CD}
V(\bar\nu)_i@>\bar{B}_{\tau}>>
\mathop{\oplus}_{\out(\tau)=i}
V(\bar\nu)_{\inn(\tau)}
@>\eps(\tau)\bar{B}_{\bar{\tau}}>>V(\bar\nu)_i\\
@VV{\bar\phi_i}V 
@V{\wr}V{}V
@VV{\bar\phi_i}V\\
V(\nu)_i@>B_{\tau}>>\mathop{\oplus}_{\out(\tau)=i}
V(\nu)_{\inn(\tau)}@>\eps(\tau)B_{{\bar\tau}}>>V(\nu)_i\\
\end{CD}
$$
Since $B\in \mu^{-1}(0)$ and $\bar B\in \mu^{-1}(0)$, 
the compositions of horizontal arrows vanish.
On the other hand $(\eps(\tau)\bb_{\bar\tau})$ is surjective 
because $\bb$ belongs to $X_0(\bar{\nu})_{i,0}$. 
Therefore, for a given $\bb$, the
fiber of $\varpi_1$ over $\bb$ is the set of elements 
$(\bar\phi,\psi,\phi'_i)$ such that the composition of
$$V(\bar\nu)_i\overset{\bar\phi_i}{\hookrightarrow} V(\nu)_i
\overset{\psi}{\rightarrow}
\Ker(\mathop{\oplus}_{\tau;\out(\tau)=i}V(\bar\nu)_{\inn(\tau)}\to 
{V(\bar\nu)_i})$$
coincides with the morphism induced by $\bar B$.
Hence the fiber ${\varpi_1}^{-1}(\bb)$ is 
connected and locally isomorphic to
$$\prod\limits_{j\in I}GL(V(\nu)_j)\times 
\Hom\Big(\C^c,
\Ker(\mathop{\oplus}_{\tau;\out(\tau)=i}V(\bar\nu)_{\inn(\tau)}\to 
{V(\bar\nu)_i})\Big)/\Hom(\C^c,V(\bar\nu)_i)\,.$$
%
\end{pf}

Now we denote by $B(\infty;\nu)$ the set of irreducible components of 
$\Lambda(\nu)$.
For $\Lambda \in B(\infty;\nu)$, we define $\eps_i(\Lambda)=\eps_i(B)$ 
by taking a generic point $B$ of $\Lambda$ . For $l\in 
{\Bbb Z}_{\geq 0}$, we set $B(\infty;\nu)_{i,\,l}$
the set of all elements of $B(\infty;\nu)$ such that $\eps_i(\Lambda)=l$.

The preceding lemma implies the following proposition.

\begin{prop}
$$B(\infty;\bar{\nu})_{i,0}\cong B(\infty;\nu)_{i,c}.$$
\end{prop}

%

\begin{defn}
Suppose that $\bar{\Lambda}\in B(\infty;\bar{\nu})_0$ corresponds to 
$\Lambda\in 
B(\infty;\nu)_c$ by this isomorphism. Then we define maps 
$\fit^c:B(\infty;\bar{\nu})_0
\to B(\infty;\nu)_c$ and $\eit^\max:B(\infty;\nu)_c\to B(\infty;\bar{\nu})_0$ 
by
$$\fit^c(\bar{\Lambda})=\Lambda,$$
$$\eit^\max(\Lambda)=\bar{\Lambda}.$$
Furthermore we define the maps 
$$\eit:\mathop{\bigsqcup}_{\nu}B(\infty;\nu)\to
\mathop{\bigsqcup}_{\nu}B(\infty;\nu)\sqcup\{0\}\text{ and}$$
$$\fit:\mathop{\bigsqcup}_{\nu}B(\infty;\nu)\to
\mathop{\bigsqcup}_{\nu}B(\infty;\nu)$$
as follows.
If $c>0$ then we define
$$\eit:B(\infty;\nu)_c @>{{\eit}^\max}>> B(\infty;\bar{\nu})_0
     @>{{\fit}^{c-1}}>> B(\infty;\nu+\alpha_i)_{c-1},$$
and  $\eit(\Lambda)=0$ for $\Lambda\in B(\infty;\nu)_0$.
We define $\fit$ by
$$\fit:B(\infty;\nu)_c @>{{\eit}^\max}>> B(\infty;\bar{\nu})_0
     @>{{\fit}^{c+1}}>> B(\infty;\nu-\alpha_i)_{c+1}.$$
\end{defn}

Then the maps $\eit^\max$ (resp. $\fit^c$) which is constructed in 
the definition may be considered as the $c$-th power of $\eit$ (resp. $\fit$).
Let us define a map $\wt:
\mathop{\bigsqcup}_{\nu}B(\infty;\nu) \to P$ by $\wt(\Lambda)=\nu \in P$ for 
$\Lambda \in B(\infty;\nu)$. We set $\varphi_i(\Lambda)=\eps_i(\Lambda)+
\langle h_i,\wt(\Lambda)\rangle$.

\begin{thm}
$\mathop{\bigsqcup}_{\nu}B(\infty;\nu)$ is a crystal in the sense of Definition
\ref{dfn:cry}.
\end{thm}

\begin{pf}
By the definition, (C1) and (C3) are automatically satisfied. 
By the definition of $\eit$, we have 
$\eps_i(\eit\Lambda)=c-1=\eps_i(\Lambda)-1$. Similarly we have
$\eps_i(\fit\Lambda)=\eps_i(\Lambda)+1$. Therefore (C2) and (C2') are 
satisfied. Since there is no $\Lambda$ such that $\eps_i(\Lambda)=
-\infty$, (C4) is satisfied.
\end{pf}

\begin{lemma}\label{ht}
If $\Lambda\in B(\infty;\nu)$ satisfies
$\varepsilon_i(\Lambda)=0$ for every $i$,
then $\nu=0$.
\end{lemma}

\begin{pf}
By the assumption,
$\mathop{\oplus}_{\inn(\tau)=i}
V(\nu)_{\out(\tau)}@>{(B_{\tau})}>>V(\nu)_i$
is surjective for a generic point $B$ of $\Lambda$.
Hence for every $n$,
$$\mathop{\oplus}_{{\sigma\in{\frak{S}}_n}\atop{\inn(\sigma)=i}}
V(\nu)_{\out(\sigma)}@>{(B_{\sigma})}>>V(\nu)_i$$
is surjective.
Then the nilpotency of $B$ implies $V(\nu)_i=0$.
\end{pf}

\subsection{}
We shall use the diagram (5.1.1) in the opposite way
to (\ref{dia:e}).
\begin{eqnarray}
&&X_0({\bar\nu})\cong X_0({-c\alpha_i})\times X_0({\bar\nu})
\mathop{\longleftarrow}^{\varpi'_1}
X'_0(-c\alpha_i,\bar\nu)
\mathop{\longrightarrow}^{\varpi'_2}
X_0(\nu).\label{dia:de}
\end{eqnarray}
We define for $B\in X_0(\nu)$
$$\eps^*_i(B)=\dimc \Ker(V(\nu)_i@>{(B_{\tau})}>>
\mathop{\oplus}_{\tau;\out(\tau)=i}
V(\nu)_{\inn(\tau)}).$$
For $\Lambda\in B(\infty;\nu)$
we define $\eps^*_i(\Lambda)$ as \(\eps_i^*(B)\)
by taking a generic point $B$ of $\Lambda$.
We set
\beqn
&&X_0(\nu)_i^p=\{B\in X_0(\nu)\,;\,\eps_i^*(B)=p\},\\
&&B(\infty;\nu)_i^p=\{\Lambda\in B(\infty;\nu)\,;\,\eit^*(\Lambda)=p\}.
\endeqn
We choose an isomorphism between $V(\nu)_i$ and its dual for every $i$.
Then $*:B\mapsto {}^tB$ gives an automorphism of
$X_0(\nu)$ and $\Lambda(\nu)$ is invariant by this automorphism.
This induces an automorphism
$*:B(\infty;\nu)\to B(\infty;\nu)$.
Since $\Lambda(\nu)$ is $G_{V(\nu)}$-invariant,
this does not depend of the choice
of isomorphisms $V(\nu)^*\simeq V(\nu)$.
The diagrams (\ref{dia:e}) and (\ref{dia:de})
are transformed by $*$.
We have
$$\eps^*_i(\Lambda)=\eps_i(\Lambda^*).$$
We define
\beqn
\eit^*{}^\max&=&*\circ \eit^\max\circ *,\\
\eit^*&=&*\circ \eit\circ *,\\
\fit^*&=&*\circ \fit\circ *,\\
\varphi_i^*(\Lambda)&=&\varphi(\Lambda^*).
\endeqn

Note that $\eit^*$ and $\fit^*$
may be defined
as $\eit$ and $\fit$
using (\ref{dia:de}) instead of (\ref{dia:e}).
We have 
\beqn
\eit^*{}^\max:B(\infty;{\nu})^c{\buildrel\sim\over\to} B(\infty;\bar\nu)^0
\endeqn

\begin{prop}
Let $\Lambda$ be an irreducible component of $\Lambda(\nu)$. We set 
$c=\eps_i^*(\Lambda)$ and $\bar\Lambda=\eit^{*}{}^\max\Lambda$. Then we have
\begin{enumerate}
\item
$$
\eps_i(\Lambda)=\max(\eps_i(\bar\Lambda),c-(\alpha_i,\bar\nu)).
$$
\item for $i\ne j$,
     $$\eps_i^*(\tilde{e_j}(\Lambda))=c,$$
     $$\eit^*{}^\max(\tilde{e_j}(\Lambda))=\tilde{e_j}(\bar\Lambda).$$
\item
     Assume $\eps_i(\Lambda)>0$. Then we have
     $$\eps_i^*(\eit(\Lambda))=
     \begin{cases}
       c &\text{if }\eps_i(\bar\Lambda)\geq c-(\alpha_i,\bar\nu),\\
       c-1 &\text{if } \eps_i(\bar\Lambda)<c-(\alpha_i,\bar\nu),
     \end{cases}
     $$
     and
     $$\eit^{*max}(\eit(\Lambda))=
     \begin{cases}
       \eit(\bar\Lambda), &\text{if }
       \eps_i(\bar\Lambda)\geq c-(\alpha_i,\bar\nu),\\
       \bar\Lambda, &\text{if } \eps_i(\bar\Lambda)<c-(\alpha_i,\bar\nu).
     \end{cases}
     $$
\end{enumerate}
\end{prop}

\begin{pf}
Since (2) may be proved in a similar way to (3)
with an easier argument,
we shall only prove (1) and (3).
Consider the diagram (\ref{dia:de}).
Let us take a generic point $\bar B$ of $\bar\Lambda$.
Then a generic point $B$ of $\varpi'_2\varpi'_1{}^{-1}(\bar B)$
is a generic point of $\Lambda$.

Fix a surjection 
$$\phi':V(\nu)\to V(\bar\nu).$$
Set
$$
N=\Coker\Big(V(\nu)_i@>{}>>V(\bar\nu)@>{(\bar{B}_{\tau})}>>
\mathop{\oplus}_{\out(\tau)=i}
V(\nu)_{\inn(\tau)}\Big).$$
Then 
$$\dimc N=\dimc\Big(\mathop{\oplus}_{\out(\tau)=i}
V(\nu)_{\inn(\tau)}\Big)-\dimc V(\bar\nu)_i
=\dimc V(\bar\nu)_i+(\alpha_i,\bar\nu).$$

The maps $B_{\bar\tau}$ induces a map
$\psi:N\to V(\bar\nu)_i$.
Let $\varphi_i:N\to V(\nu)_i$ be a generic map
such that
$$
\psi=\phi'_i\circ \varphi.
$$
Then $B$ is given as follows.
\beqn
B_\tau&=&
\begin{cases}
\bar B{\,}_\tau&\text{if $\inn(\tau),\out(\tau)\not=i$,}\\
V(\nu)_i@>{\phi'_i}>>
V(\bar\nu)_i@>{\bar B_\tau}>>V(\bar\nu)_{\inn(\tau)}
&\text{if $\out(\tau)=i$,}\\
V(\nu)_{\inn(\bar\tau)}@>{}>>
\mathop{\oplus}\limits_{\tau';\,\out(\tau')=i}
V(\nu)_{\inn(\tau')}@>{}>>N@>{\varphi}>>V(\nu)_i
&\text{if $\inn(\tau)=i$.}
\end{cases}
\endeqn

Since $\varphi$ is generic, we have
\beq
\dimc\Ker\varphi&=&\max(\dimc \Ker \psi-c,0)
\endeq

Let us calculate
$\eps_i(B)$.
Since
$\Im\Big(\mathop{\oplus}_{\inn(\tau)=i}
V(\nu)_{\out(\tau)}@>{}>>V(\nu)_i\Big)
=\Im(N@>{\varphi}>>V(\nu)_i)$,
\beqn
\eps_i(B)&=&\dimc V(\nu)_i-\dimc \Im(\varphi)\\
&=&\dimc V(\nu)_i-\dimc N+\dim \Ker(\varphi)\\
&=&\max\Big(\dim V(\nu)_i-\dimc N+\dim\Ker(\psi)-c\,
,\,\dim V(\nu)_i-\dim N\Big)\\
&=&\max\Big(\dim\Coker \psi,\,c-(\alpha_i,\bar\nu)\Big)\\
&=&\max\Big(\eps_i(\bar\Lambda),\,c-(\alpha_i,\bar\nu)\Big)
\endeqn
Thus we obtain (1).

Now let us prove (3).
Set $\nu'=\nu+\alpha_i$ and $\bar\nu'=\bar\nu+\alpha_i$.
Let us take a generic hyperplane
$H$ of $V(\nu)_i$ containing $\Im(\varphi)$.
Then taking $V(\nu')\simeq H$,
we obtain a generic point $B'$ of $\eit(\Lambda)$.
Similarly, taking a generic hyperplane $\bar H$ of $V(\bar\nu)_i$
containing $\Im (\psi)$,
we obtain a generic point of
$\eit(\bar\Lambda)$.

If $\dim\Ker \psi\ge c$ (i.e. $\eps_i(\bar\Lambda)\ge c-(\alpha_i,\bar\nu)$),
then
$\Im(\varphi)\supset \Ker(\phi'_i)$
and hence
$\phi'_i(H)$ is a hyperplane of $V(\bar\nu)_i$.
Therefore
$\eps_i^*(\te_i(\Lambda))=c$
and
$\eit^*{}^\max(\eit(\Lambda))=\eit(\bar\Lambda)$.

If $\dim\Ker \psi<c$, then
$\Im(\varphi)\not\supset \Ker(\phi'_i)$
and hence
$\phi'_i(H)=V(\bar\nu)_i$.
Therefore
$\eps_i^*(\te_i(\Lambda))=c-1$
and
$\eit^*{}^\max\eit(\Lambda)=\bar\Lambda$.
\end{pf}

We recall that $B(\infty)$ is the crystal base of $U_q^-({\frak g})$.

\begin{thm}
We have an isomorphism of crystals
$$\mathop{\bigsqcup}_{\nu\in Q_-}B(\infty;\nu)\cong B(\infty)\text{ }.$$
\end{thm}

\begin{pf}
We define a map $\Phi_i:\mathop{\bigsqcup}_{\nu}B(\infty;\nu)\to
\mathop{\bigsqcup}_{\nu}B(\infty;\nu)\otimes B_i$ by $\Lambda \mapsto
\eit^\max(\Lambda)\otimes\fit^{\eps^*_i(\Lambda)}b_i$. 
It is clear that this map is well-defined 
and injective. Moreover,
it is a strict morphism of crystals by the preceding lemma.
Now we can apply Proposition \ref{char}
because the condition (\ref{cond:1}) is satisfied by
Lemma \ref{ht}.
\end{pf}

We denote by $\Lambda_b\in\mathop{\bigsqcup}_{\nu\in Q_-}B(\infty;\nu) $ 
the corresponding element to $b\in B(\infty)$ under this isomorphism. 
The following proposition is proved by Lusztig.
\begin{prop}
$\Lambda(\nu)$ is a Lagrangian subvariety of $X_0(\nu)$.
\end{prop}

By this result, any $\Lambda_b$ in $B(\infty;\nu)$
is an irreducible Lagrangian subvariety of $X_0(\nu)$.

\section{Review of the theory of canonical base}
\subsection{Canonical base}
Let us recall the results on Lusztig on canonical bases.
We write $\D(X)$ for the bounded derived category of complexes of
sheaves of $\C$-vector spaces on the associated complex variety
with an algebraic variety $X$ over $\C$. Objects of
$\D(X)$ are referred to as complexes. We shall use the notations of [BBD];
in particular, $[d]$ denotes a shift by $[d]$ degrees,
and for a morphism $f$ of algebraic varieties, $f^*$ denotes the
inverse image functor, $f_!$ denotes direct image with compact support,
etc.

We fix an orientation  $\om$ of quiver.
Let $\nu\in Q_-$ and let $S_{\nu}$ be the set of all 
pairs $(\text{i, a})$ where $\text{i}=
(i_1,i_2,\cdots,i_m)$ is a sequence of elements of $I$ and 
$\text{a}=(a_1,a_2,\cdots,a_m)$ is a sequence of non-negative 
integers such that $\mathop{\sum}_{l}a_l\alpha_{i_l}=-\nu$. 
Now let $V\in {\cal V}_\nu$ and let $(\text{i},\text{a})\in
S_{\nu}$. A flag of type $(\text{i},\text{a})$ is, 
by definition, a sequence $\phi=(V=V^0\supset V^1\supset\cdots\supset 
V^m=0)$ of $I$-graded
subspace of $V$ such that, for any $l=1,2,\cdots,m$, the $I$-graded vector
space $V^{l-1}/V^l$ is zero in degrees $\ne i_l$ and has dimension $a_l$ in
degree $i_l$. We define a variety $\tilde{\cal F}_{\text{i},
\text{a}}$ of all pairs $(B,\phi)$ such that $B\in \ev$ and 
$\phi$ is a $B$-stable flag of type $(\text{i},
\text{a})$. The group $G_V$ acts on 
$\tilde{\cal F}_{\text{i},\text{a}}$ in natural way.
We denote by $\pi_{\text{i},\text{a}}:
\tilde{\cal F}_{\text{i},\text{a}}\to \ev$ the 
natural projection. We note that $\pi_{\text{i},\text{a}}$
is a $G_V$-equivariant proper morphism. We set $L_{\text{i},
\text{a};\om}=(\pi_{\text{i},\text{a}})_!
(1)\in \D(\ev)$.
Here $1\in \D(\tilde{\cal F}_{\text{i},\text{a}})$ 
is the constant sheaf on 
$\tilde{\cal F}_{\text{i},\text{a}}$.
By the decomposition theorem [BBD], $L_{\text{i},\text{a}
;\om}$ is a semisimple complex.

Let $\pv$ be the set of isomorphism class of simple perverse sheaves $L$ on
$\ev$ such that $L[d]$ appears as direct summand of 
$L_{\text{i},\text{a};\om}$ for some 
$(\text{i},\text{a})\in S_{\nu}$ and
some $d\in {\Bbb Z}$. We write $\qv$ for the subcategory of $\D(\ev)$
consisting of all complexes that are isomorphic to finite direct sums of
complexes of the form $L[d]$ for various simple perverse sheaves $L\in \pv$
and various $d\in {\Bbb Z}$. Any complex in $\qv$ is semisimple and
$G_V$-equivariant.

Take $V\in {\cal V}_{\nu},V'\in {\cal V}_{\nu'}
,\bar{V}\in {\cal V}_{\bar{\nu}}$ for ${\nu}={\nu'}
+{\bar{\nu}}$ in $Q_-$). We consider the diagram
\beq
&&E_{\bar V,\om}\times E_{V',\om}\mathop{\longleftarrow}^{p_1}E'
\mathop{\longrightarrow}^{p_2}E''\mathop{\longrightarrow}^{p_3}\ev\,.
\label{base}
\endeq
Here $E'$ is the variety of
$(B,\bar{\phi},\phi')$ where $B\in E_{V,\Omega}$ 
and $0\to \bar V\buildrel \bar\phi\over \longrightarrow
V\buildrel\phi'\over\longrightarrow V'\to 0$
is a B-stable exact sequence of $I$-graded vector spaces,
and $E''$ is the variety of $(B,C)$
where $B\in E_{V,\Omega}$ and $C$ is a $B$-stable $I$-graded subspace of
$V$ with $\dim C=\bar\nu$.
%
The morphisms $p_1$, $p_2$ and $p_3$ are defined by
$p_1(B,\bar\phi,\phi')=(B|{}_{\bar V},B|_{V'})$,
$p_2(B,\bar\phi,\phi')=(B,\Im(\bar\phi))$ and $p_3(B,C)=B$.
Note that
$p_1$ is smooth with connected fiber, $p_2$ is a principal $G_{V'}\times
G_{\bar{V}}$-bundle, and $p_3$ is proper.

Let $L'\in {\cal Q}_{V',\om}$ and $\bar{L}\in {\cal Q}_{\bar{V},\om}$. 
Consider the exterior tenser product $\bar L\boxtimes L'$.
Then there is $(p_2)_{\flat}p_1^*(\bar L\boxtimes L')\in\D(E'')$
such that  $(p_2)^*(p_2)_{\flat}p_1^*(\bar L\boxtimes L')
\cong p_1^*(\bar L\boxtimes L')$.
We define $L'*\bar{L}\in \qv$ by 
$(p_3)_!(p_2)_{\flat}p_1^*(\bar L\boxtimes L')[d_1-d_2]$ 
where $d_i$ is the fiber dimension of $p_i$ ($i=1,2$).
Let ${\cal K}_{V,\om}$ be the Grothendieck group of $\qv$.
We considered as a ${\Bbb Z}[q,q^{-1}]$-module by 
$q(L)=L[1]$, $q^{-1}(L)=L[-1]$. 
Then ${\cal K}_{\om}=\mathop{\oplus}_{\nu\in Q_-}{\cal K}_{V(\nu),\Omega}$
has a structure of an associative graded ${\Bbb Z}[q,q^{-1}]$-algebra
by the operation $*$.
We denote by $F_i\in{\cal K}_{V(-\alpha_i),\om}$ 
the element attached to $1\in \D(E_{V(-\alpha_i),\om})$.

\begin{thm}\cite{L;canI}
There is a unique ${\Bbb Q}(q)$-algebra isomorphism 
$$\Gamma_{\om}:U_q^-({\frak g})\to {\cal K}_{\om}
\mathop{\otimes}_{{\Bbb Z}[q,q^{-1}]}{\Bbb Q}(q)$$
such that $\Gamma_{\om}(f_i)=F_i$.
\end{thm}

Let us identify $L\in \pv$ with $L\otimes 1\in {\cal K}_{\om}
\mathop{\otimes}_{{\Bbb Z}[q,q^{-1}]}{\Bbb Q}(q)$. We set $\text{\bf B}
=\Gamma_{\om}^{-1}(\bigsqcup_{V\in {\cal V}}\pv)$ and call it
the canonical basis of $U_q^-({\frak g})$. 
By [GL], $\text{\bf B}$ and $B(\infty)$ are canonically identified.
For $b\in B(\infty)$
the corresponding perverse sheaf is denoted by
$L_{b,\om}$.

\subsection{}
Let $Y$ be a smooth algebraic variety. For any $L\in \D(Y)$,
we denote by $SS(L)$ the singular support 
(or the characteristic variety) of $L$.
It is known that $SS(L)$ is a closed Lagrangian subvariety of $T^*Y$
(See \cite{KS}).

We recall that $T^*\ev$ is identified with $X_V$.
By the Fourier transform method, we have
\begin{thm}\label{inv}\cite{L;quiver}
$SS(L_{b,\om})$ does not depend on the choice of $\om$.
\end{thm}

We say $i\in I$ is sink (resp. source) of $\om$ if there is no arrow
$i\to j$ (resp. $j\to i$) in $\om$.
\begin{thm}\label{estch}
\begin{enumerate}
\item For any $L\in \pv$ the singular support $SS(L)$ is a union of 
     irreducible components of $\lv$.
\item For any $b\in B(\infty)$ and $i\in I$, we have
\end{enumerate}
\beq
&&\Lambda_b\subset SS(L_{b,\om})\subset \Lambda_b\cup 
\mathop{\bigcup}_{\eps_i(b')>\eps_i(b)}\Lambda_{b'}.
\label{char.est}
\endeq
\end{thm}
\begin{pf}
The first statement is due to Lusztig \cite{L;quiver}.
By taking $\om$ such that $i$ is a sink,
the second statement follows from
\cite{L;canI}.
\end{pf}

Note that if there is a bijection $s:B(\infty)\to B(\infty)$ such that
$SS(L_{b,\om})\supset \Lambda_{s(b)}$ for any $b\in B(\infty)$, then $s$
must be the identity ({\it cf.} Problem in \cite{L;quiver}).
In fact, by the decreasing induction on $\eps_i(b)$, 
(\ref{char.est}) implies $s(b)=b$. 

The following problem is also asked by Lusztig \cite{L;quiver}.

\begin{prob}
If the underlying graph is of type $A,D,E$, then the singular support of any
$L\in \pv$ is irreducible.
\end{prob} 

Furthermore he noted that the next conjecture \cite{KL;sp}
follows from Problem 1 for type $A$ (see \S\ref{equiv}).
In fact it is easy to see that they are equivalent.

\begin{conj} Let $X$ be the flag manifold for $SL(n)$
and let $X_w$ be the Schubert variety of $SL(n)$ which corresponds to the 
element $w$ of the Weyl group $W$. Then the singular support of 
$\sideset{^{\pi}}{_{X_w}}{\Bbb C}$ is irreducible.  
\end{conj}

In the next section we construct a counterexample of Problem 1 for a 
graph of type $A$. 

%
%

\section{Counterexample to problem 1}
\subsection{}
In this and the next section we assume 
that the underlying graph is of type $A$.

Let us take $\nu\in Q_-$ and $V\in {\cal V}_{\nu}$.
Let ${\cal O}_{\om}$ be a $G_V$-orbit in $\ev$. 
As the underlying graph is of type $A$,
$\ev$ has finitely many $G_V$-orbits. 
By \cite{L;canI} we know that there is one-to-one correspondence between
$G_V$-orbits ${\cal O}$ in $\ev$ 
between the crystal basis $b\in B(\infty)$ of $U_q^-({\frak g})$ of weight 
$\nu$ by $\Lambda_b=T^*_{\cal O}\ev$.
For $b\in B(\infty)$, we denote by $\CO_{b,\om}$ the corresponding
$G_V$-orbit.
The next theorem is due to Lusztig (See \cite{L;canI}.).

%
\begin{thm}
Let $b\in B(\infty)$. Then we have $$L_{b.\om}=
       \sideset{^{\pi}}{_{{\cal O}_{b,\om}}}{\Bbb C}$$
where ${\Bbb C}_{{\cal O}_{b,\om}}$ is the constant sheaf on 
${\cal O}_{b,\om}$
and $\sideset{^{\pi}}{}{\cdot}$ is the minimal extension functor.
\end{thm}

Note that $SS(L_{b.\om})$ depends only on $b\in B(\infty)$ and not on 
$\om$ (cf. Theorem {\ref{inv}}).
\subsection{}
In the rest of the section, we shall present a counterexample 
of Problem 1 when the underlying graph is of type $A_5$. 
Let us take a graph  of type $A_5$ and its orientation $\om$
as follows; 
$$\mathop{\circ}^1\mathop{\longleftarrow}^{\tau_1}
\mathop{\circ}^2\mathop{\longleftarrow}^{\tau_2}
\mathop{\circ}^3\mathop{\longleftarrow}^{\tau_3}
\mathop{\circ}^4\mathop{\longleftarrow}^{\tau_4}
\mathop{\circ}^5.$$

Let $\nu=-2\alpha_1-4\alpha_2-4\alpha_3-4\alpha_4-2\alpha_5$.
Set $b=\ft_2\ft_1\ft_3\ft_2\ft_4^2\ft_3^2\ft_2\ft_1\ft_5^2\ft_4^2\ft_3\ft_2
u_{\infty}$ 
and $b'=\ft_2^2\ft_1^2\ft_3^2\ft_4^2\ft_3^2\ft_2^2\ft_5^2\ft_4^2u_{\infty}$.
Then the following points $B_0$ and $B'_0$ of $\ev$
are in $\CO_{b,\om}$ and $\CO_{b',\om}$, respectively.
\beqn
&&(B_0)_{\tau_1}=
\begin{pmatrix}
1 & 0 & 0 & 0 \\
0 & 0 & 1 & 0
\end{pmatrix},
\text (B_0)_{\tau_2}=
\begin{pmatrix}
1 & 0 & 0 & 0 \\
0 & 1 & 0 & 0 \\
0 & 0 & 0 & 0 \\
0 & 0 & 0 & 1
\end{pmatrix},\\
&&\phantom{**********}
(B_0)_{\tau_3}=
\begin{pmatrix}
1 & 0 & 0 & 0 \\
0 & 1 & 0 & 0 \\
0 & 0 & 1 & 0 \\
0 & 0 & 0 & 0
\end{pmatrix}
\text{, } (B_0)_{\tau_4}=
\begin{pmatrix}
0 & 0 \\
1 & 0 \\
0 & 0 \\
0 & 1
\end{pmatrix},\\
&&(B'_0)_{\tau_1}=
\begin{pmatrix}
0 & 0 & 1 & 0 \\
0 & 0 & 0 & 1
\end{pmatrix},
\text (B'_0)_{\tau_2}=
\begin{pmatrix}
1 & 0 & 0 & 0 \\
0 & 1 & 0 & 0 \\
0 & 0 & 0 & 0 \\
0 & 0 & 0 & 0
\end{pmatrix},\\
&&\phantom{**********}
(B'_0)_{\tau_3}=
\begin{pmatrix}
1 & 0 & 0 & 0 \\
0 & 1 & 0 & 0 \\
0 & 0 & 0 & 0 \\
0 & 0 & 0 & 0
\end{pmatrix}
\text{, } (B'_0)_{\tau_4}=
\begin{pmatrix}
0 & 0 \\
0 & 0 \\
1 & 0 \\
0 & 1
\end{pmatrix}.
\endeqn

Now we can state a counterexample of Conjecture 1.
\begin{thm} \label{counter}
$$SS(\sideset{^{\pi}}{_{{\cal O}_{b,\om}}}{\Bbb C})\supset \Lambda_b\cup
\Lambda_{b'}.$$
\end{thm}

\begin{rem}
In fact, although we don't give a proof
(relying on Lemmas \ref{crt1} and \ref{crt2}), they coincide.
\end{rem}

Let $U$ be the subvariety of $\ev$ consisting of
elements $B=(B_{\tau_1},B_{\tau_2},B_{\tau_3},B_{\tau_4})\in \ev$
of the form
$B_{\tau_1}=
\begin{pmatrix}
X_1, & I_2 
\end{pmatrix}
$, 
$B_{\tau_2}=
\begin{pmatrix}
I_2, & Y_1 \\
0,   & Y_2 
\end{pmatrix}
$, 
$B_{\tau_3}=
\begin{pmatrix}
I_2 ,& Z_1 \\
0 ,& Z_2
\end{pmatrix}
$ and 
$B_{\tau_4}=
\begin{pmatrix}
0 \\
I_2
\end{pmatrix}
$
where $I_2=\begin{pmatrix}1&0\\0&1\end{pmatrix}$ and $X_1$, $Y_1$, $Y_2$,
$Z_1$, $Z_2$ are $2\times 2$-matrices.

Let $M$ be the complex vector space of 
$(A_i)_{i\in {\Bbb Z}/4{\Bbb Z}}$ 
where each $A_i$ is a $2\times2$-matrix, let $\overline{S}$ be the 
closed subvariety of $M$ consisting all elements such that $A_{i+1}A_i=0$ 
and $\text{rank}(A_i)\leq 1$ and let $S$ be the subvariety of $M$ consisting 
of elements such that $A_{i+1}A_i=0$ and $\text{rank}(A_i)=1$. It is clear 
that $\overline{S}$ is the closure of $S$.

\begin{lemma}
\begin{enumerate}
\item ${\cal O}_{b,\om}\cap U\cong S,$
\item ${\cal O}_{b',\om}\cap U=\{B'_0\}$ and ${\cal O}_{b',\om}$ and 
     $U$ intersect transversally.
\end{enumerate}
\end{lemma}

\begin{pf}
Let us prove (1). $B=(B_{\tau_1},B_{\tau_2},B_{\tau_3},B_{\tau_4})\in
{\cal O}_{b,\om}\cap U$ is characterized following properties;
\renewcommand{\theenumi}{\roman{enumi}}
\renewcommand{\labelenumi}{(\theenumi)}
\begin{enumerate}
\item $\rank B_{\tau_1}=\rank(X_1,\,I_2)=2$,
\item $\rank B_{\tau_2}B_{\tau_1}=\rank(X_1\,,X_1Y_1+Y_2)=1$,
\item $\rank B_{\tau_3}B_{\tau_2}B_{\tau_1}=
      \rank(X_1,\,X_1Z_1+(X_1Y_1+Y_2)Z_2)=1$,
\item $\rank B_{\tau_4}B_{\tau_3}B_{\tau_2}B_{\tau_1}=
      \rank\big(X_1Z_1+(X_1Y_1+Y_2)Z_2\big)=0$,
\item $\rank B_{\tau_2}=\rank 
      \begin{pmatrix}
       I_2, & Y_1 \\
       0   ,& Y_2
      \end{pmatrix}
      =3$,
\item $\rank B_{\tau_3}B_{\tau_2}=\rank
      \begin{pmatrix}
       I_2, & Z_1+Y_1Z_2 \\
       0   ,& Y_2Z_2
      \end{pmatrix}
      =2$,
\item $\rank B_{\tau_4}B_{\tau_3}B_{\tau_2}=\rank
      \begin{pmatrix}
       Z_1+Y_1Z_2 \\
       Y_2Z_2
      \end{pmatrix}
      =1$,
\item $\rank B_{\tau_3}=\rank
      \begin{pmatrix}
       I_2, & Z_1 \\
       0   ,& Z_2
      \end{pmatrix}
      =3$,
\item $\rank B_{\tau_4}B_{\tau_3}=\rank
      \begin{pmatrix}
       Z_1 \\
       Z_2
      \end{pmatrix}
      =1$,
\item $\rank B_{\tau_4}=
      \begin{pmatrix}
       0\\
       I_2 
      \end{pmatrix}
      =2$.
\end{enumerate}

Set $\tilde{Z}_1=Z_1+Y_1Z_2$. By a direct calculation, we conclude that 
${\cal O}_{b,\om}\cap U$ is  isomorphic to the set of 
$(X_1,Y_2,\tilde{Z}_1,Z_2)$ such that 
$$\rank X_1=\rank Y_2=\rank \tilde{Z}_1=\rank Z_2=1,$$
$$(\cof Y_2)X_1=X_1\tilde{Z}_1=Y_2Z_2=\tilde{Z}_1(\cof Z_2)=0.$$
Here we denote the cofactor of $A$ by $\cof A$. This proves (1).

The similar arguments yield (2).
\end{pf}


To see Theorem 7.2.1, it is enough to show that
$$SS(\sideset{^{\pi}}{_{S}}{\Bbb C})\text{ {\it is not 
an irreducible variety.}}\leqno(7.2.1)$$

Let $\sideset{}{_{S|M-\partial{S}}}{\cal B}$ be the $D$-module of the delta 
function on $S$ in $M-\partial{S}$ where $\partial{S}=\overline{S}-S$.
By the Riemann-Hilbert correspondence, (7.2.1) is equivalent to
$$SS(\sideset{^{\pi}}{_{S|M-\partial{S}}}{\cal B})\text{ {\it is not
an irreducible variety.}}
\leqno(7.2.2)$$

\subsection{}
In this subsection we shall prove (7.2.2).

Assuming that 
$SS(\sideset{^{\pi}}{_{S|M-\partial{S}}}{\cal B})=\overline{T^*_{S}M}$,
we shall deduce a contradiction.

We denote by $\sideset{^{F}}{}{\cdot}\,$ the Fourier transformation functor.
Then there is an isomorphism $F:\Gamma(M;{\frak M})\to\Gamma(M^*;
\sideset{^{F}}{}{\frak M})$ such that $F\circ x_i=\partial_{\xi_i}\circ F$,
etc.. Here $\frak M$ is a $D$-module on $M$, $x_i$ is a local coordinate of 
$M$ and $\xi_i$ is the corresponding dual coordinate.

According to [KS] we have 
$$SS(\sideset{^{\pi}}{_{S|M-\partial{S}}}{\cal B})
=SS(\sideset{^{F\pi}}{_{S|M-\partial{S}}}{\cal B})$$
under the identification $T^*M\cong T^*M^*$  where $M^*$ is the dual space 
of $M$. On the other hand, it is known that
$$\text{supp}(\sideset{^{F\pi}}{_{S|M-\partial{S}}}{\cal B})
=SS(\sideset{^{F\pi}}{_{S|M-\partial{S}}}{\cal B})\cap T^*_{M^*}M^*.$$ 
Set $S^*=\text{supp}(\sideset{^{F\pi}}{_{S|M-\partial{S}}}{\cal B}).$
Then $S^*$ is the polar variety of $S$, i.e.
$$S^*=\overline{T^*_{S}M}\cap T^*_{M^*}M^*\subset M^*.$$


Let $(A_i)\in M$ and let 
$A_i=
\begin{pmatrix}
x_i & y_i \\
z_i & w_i
\end{pmatrix}$
for $i\in {\Bbb Z}/4{\Bbb Z}$. We note that $x_i$, $y_i$, $z_i$, $w_i$
are coordinates of $M$. Let $\xi_i$, $\eta_i$, $\zeta_i$, $\rho_i$ 
be the corresponding dual coordinates of $M^*$. 
Set $A_i^*=
\begin{pmatrix}
\xi_i  & \zeta_i \\
\eta_i & \rho_i
\end{pmatrix}.$

By an easy calculation, we see that $S^*$ is the variety of  
$(A_i^*)_{i\in {\Bbb Z}/4{\Bbb Z}}\in M^*$ such that two eigenvalues of the 
$2\times 2$-matrix $A_1^*A_2^*A_3^*A_4^*$ coincide.
That is, 
$$S^*=\{(A_i^*)\in M^*\,;\,f=0\}$$
where 
$f=(\theta_{1,1}-\theta_{2,2})^2
+4\theta_{1,2}\theta_{2,1}$ with
$\begin{pmatrix}
\theta_{1,1} & \theta_{1,2} \\
\theta_{2,1} & \theta_{2,2}
\end{pmatrix}
=A_1^*A_{2}^*A_{3}^*A_{4}^*$.
\hb
Hence $SS(\sideset{^{F\pi}}{_{S|M-\partial{S}}}{\cal B})$ is an 
irreducible $D_{M^*}$-module supported on $S^*$.
There is a $({\Bbb C}^*)^4\times \text{GL}
({\Bbb C}^2)^4$-action on $M$, by which
$((c_i),(g_i))\in({\Bbb C}^*)^4 \times \text{GL}({\Bbb C}^2)^4$ sends
$(A_i)$ to $(c_ig_{i+1}A_ig_i^{-1})$.
Then $S$ is invariant by this action and hence 
$\sideset{^{\pi}}{_{S|M-\partial{S}}}{\cal B}$ is $({\Bbb C}^*)^4\times
\text{GL}({\Bbb C}^2)^4$-equivariant. Hence its Fourier transform
$\sideset{^{F\pi}}{_{S|M-\partial{S}}}{\cal B}$ is also 
$({\Bbb C}^*)^4\times\text{GL}({\Bbb C}^2)^4$-equivariant.

%
%
Let $S^*_0$ be a unique open 
$({\Bbb C}^*)^4\times\text{GL}({\Bbb C}^2)^4$-orbit of $S^*$.
Then its isotropy subgroup is connected.
Since $\text{supp}(\sideset{^{F\pi}}{_{S|M-\partial{S}}}{\cal B})=S^*$ and
$\sideset{^{F\pi}}{_{S|M-\partial{S}}}{\cal B}$ is irreducible, we have 
$$\sideset{^{F\pi}}{_{S|M-\partial{S}}}{\cal B}=
\sideset{^{\pi}}{}({\cal B}_{S^*_0|M^*-\partial{S^*_0}}\otimes{L})$$
where $L$ is an irreducible 
$({\Bbb C}^*)^4\times\text{GL}({\Bbb C}^2)^4$-equivariant local system on 
$S^*_0$. As the isotropy subgroup of $S^*_0$ is connected, any irreducible
$({\Bbb C}^*)^4\times\text{GL}({\Bbb C}^2)^4$-equivariant local system on
$S^*_0$ must be trivial. Therefore we have 
$\sideset{^{F\pi}}{_{S|M-\partial{S}}}{\cal B}
=\sideset{^{\pi}}{_{S^*_0|M^*-\partial{S^*_0}}}{\cal B}$.

The next result is due to Barlet-Kashiwara [BK].
\begin{prop}[Barlet-Kashiwara] 
\begin{enumerate}
\item $\sideset{^{\pi}}{_{S^*_0|M^*-\partial{S^*_0}}}{\cal B}\subset 
     {\cal O}_{M^*}[1/f]/{\cal O}_{M^*}$.
\item $\frac{\partial f}{\partial \xi_i}\delta(f)\in 
     \sideset{^{\pi}}{_{S^*_0|M^*-\partial{S^*_0}}}{\cal B}$ where
     $\delta(f)=1/f\ \mod{\cal O}_{M^*}$ is the delta 
     function.
\end{enumerate}
\end{prop}

Since Proposition 7.3.1, $\text{supp}F^{-1}(\partial_{\xi_i}f)\delta(f)\subset
\bar{S}$ and
$\text{det}A_i$ vanishes on $\bar{S}$, Hilbert-Nullstellensatz 
guarantees the existence of a positive integer $m$ such that
$$
\begin{vmatrix}
{x_i}  & {y_i} \\
{z_i} & {w_i}
\end{vmatrix}^m
F^{-1}((\partial_{\xi_i} f)\delta(f))=0.
$$
Applying the Fourier transformation, we have 
$$
\begin{vmatrix}
\partial_{\xi_i}  & \partial_{\zeta_i} \\
\partial_{\eta_i} & \partial_{\rho_i}
\end{vmatrix}^m
(\partial_{\xi_i} f)\delta(f)=0.\leqno(7.3.1)
$$
On the other hand, a direct calculation leads
$$
\begin{vmatrix}
\partial_{\xi_i}  & \partial_{\zeta_i} \\
\partial_{\eta_i} & \partial_{\rho_i}
\end{vmatrix}
(\partial_{\xi_i} f)\delta^{(k)}(f)=(4k-2)d(\partial_{\xi_i} f)
\delta^{(k+1)}(f),
$$
where $d=\text{det}(A_{i+1}^*A_{i+2}^*A_{i+3}^*)$.
As $d$ is a polynomial free from $\xi_i$, $\zeta_i$, $\eta_i$,
$\rho_i$, we get
$$
\begin{vmatrix}
\partial_{\xi_i}  & \partial_{\zeta_i} \\
\partial_{\eta_i} & \partial_{\rho_i}
\end{vmatrix}^m
(\partial_{\xi_i} f)\delta(f)
=\Big(\mathop{\prod}_{k=0}^{m-1}(4k-2)\Big)d^m(\partial_{\xi_i} f)
\delta^{(m)}(f).
$$
Since $k$ is a integer, the right hand side never vanishes. This contradicts
(7.3.1).
Thus we complete the proof of Theorem \ref{counter}.
\section{Relation with Schubert cells}
\subsection{ }\label{equiv}
We consider the Dynkin diagram of type $A_{2n-1}$ and take its orientation
$\Omega_0$ as follows:
$$\om_0;\text{ }\mathop{\circ}^1@<\tau_1<< \mathop{\circ}^2 @<\tau_2<<
\cdots@<\tau_{2n-2}<<\mathop{\circ}^{2n-1}.$$

Let $\nu_{cl}=\sum_{i=1}^{2n-1}-\nu_{cl}(i)\alpha_i$ where $\nu_{cl}(i)=i$ 
(for $1\leq i\leq n$), $=2n-i$ (for $n\leq i\leq 2n-1$) and let  $V_{cl}\in 
{\cal V}_{\nu_{cl}}$.

Let us set
$$E_{V_{cl}}^{\natural}=\{B\in E_{V_{cl},\om_0}\,;\,B_{\tau_i}
\text{ is injective
for }1\leq i\leq n-1 \text{ and surjective } n\leq i\leq 2n-2\}.$$
It is clear that $E_{V_{cl}}^{\natural}$
is $G_V$-invariant. 

Let $G$ be $GL(n,{\Bbb C})$, $B$ a Borel subgroup of $G$, $W$ the Weyl group
of $G$ and $X=G/B$ the flag variety. We set $X_w=BwB/B$ ($w\in W$).
Then $X=\bigsqcup_{w\in W}X_w$ gives a cellular 
decomposition of $X$.

The decomposition of $X\times X$ to $G$-orbits is given by $X\times X
=\bigsqcup_{w\in W}Y_w$ with $Y_w=G\cdot (\{eB\}\times X_w)$. Then, 
the following two conditions are equivalent:
$$SS(\sideset{^{\pi}}{_{X_w}}{\Bbb C})\text{ {\it is an irreducible variety}.}
\leqno(8.1.1)$$
$$SS(\sideset{^{\pi}}{_{Y_w}}{\Bbb C})\text{ {\it is an irreducible variety}.}
\leqno(8.1.2)$$

We have a $G$-equivariant isomorphism
$$E_{V_{cl}}^{\natural}/\prod_{j\not=n}GL(V_{cl}{}_j)
\simeq X\times X.$$
Therefore there is a one-to-one correspondence between
$G$-orbits of $X\times X$ and
$G_{V_{cl}}$-orbits of $E_{V_{cl}}^{\natural}$.
Let us denote by $\CO_{w,\om_0}$ the
$G_{V_{cl}}$-orbit of $E_{V_{cl}}^{\natural}$
corresponding to $Y_w$.
Then we have
\beq
\text{
The irreducibility of 
$SS(\sideset{^{\pi}}{_{X_w}}{\Bbb C})$ is equivalent to
that of $SS(\sideset{^{\pi}}{_{{\cal O}_{w,\om_0}}}{\Bbb C})$.}
\endeq

\subsection{}
For an orientation $\om$ we say 
that $i\in I$ is sink (resp. source) of $\om$ 
if there is no arrow
$i\to j$ (resp. $j\to i$) in $\om$. 
\begin{lemma}\label{crt1}
\begin{enumerate}
\item Let $b\in B(\infty)$. 
     If $SS(L_{b,\om}) \supset
     \Lambda_{b'}$, then $\eps_i(b)\leq \eps_i(b')$ for any $i\in I$.
\item If $\eps_i(b)=\eps_i(b')$, then the condition 
     $SS(L_{b,\om}) \supset
     \Lambda_{b'}$ is equivalent to
     $$SS(L_{{\tilde{e}_i}^\max b,\om})
     \supset\Lambda_{{\tilde{e}_i}^\max b'}.$$
\end{enumerate}
\end{lemma}

\begin{pf}
(1) is already seen in Theorem \ref{estch}.
Let us prove (2).
Let us choose an orientation
$\om$ such that $i$ is a sink.
Set $\nu=\wt(b)$, $m=\eps_i(b)$ and $\bar\nu=\nu+m\alpha_i$.
Let $Z$ be the subvariety of $E_{V(\nu),\om}$ consisting of $B$ such that
$$\mathop\oplus\limits_{{\tau\in \om}\atop{\inn(\tau)=i}}V(\nu)_{\out(\tau)}\to V(\nu)_i$$
has cokernel of dimension $m$.
Similarly let $Y$ be the subvariety of $E_{V(\bar\nu),\om}$ 
consisting of $\bar B$ such that
$$\mathop
\oplus\limits_{{\tau\in \om}\atop{\inn(\tau)=i}}V(\bar\nu)_{\out(\tau)}
\to V(\bar\nu)_i$$
is surjective.

Then there is a $GL(V(\bar\nu)_i)$-bundle
$p:Z\to Y$.
Then $\CO_{\eit^mb,\om}$ and $\CO_{\eit^mb',\om}$
are contained in $Y$,
and
$$\CO_{b,\om}=p^{-1}(\CO_{\eit^\max b,\om})\ \text{ and }\ 
\CO_{b',\om}=p^{-1}(\CO_{\eit^\max b,\om}).$$
This shows immediately (2).
\end{pf}

For an orientation $\om$, let $s_i\om$ ($i\in I$) be the 
orientation obtained from $\om$ by reversing each arrow that ends 
or starts at $i$.

We define a map $S_i:\{b\in B(\infty)\,;\,\eps_i(b)=0\}\to
\{b\in B(\infty)\,;\,\eps_i^*(b)=0\}$ by $S_i(b)=
\tilde{f}_i^{\varphi_i^*(b)}\tilde{e}_i^*{}^\max b$.
Then $S_i$ is bijective. Note that
$\wt(S_i(b))=s_i(\wt(b))$ (see [S]). Here $s_i$ is the simple reflection.

\begin{lemma}\label{crt2}
Assume that $b,b'\in B(\infty)$ has the same weight and
that $i\in I$ satisfies $\eps_i(b)=\eps_i(b')=0$.
Then the following
two conditions are equivalent;
\begin{enumerate}
\item $SS(L_{b,\om}) \supset \Lambda_{b'}$,
\item $SS(L_{S_i(b),s_i\om}) \supset\Lambda_{S_i(b')}$.
\end{enumerate}
\end{lemma}

\begin{pf}
We choose an orientation
$\om$ such that $i$ is a sink.
Set $\nu=\wt(b)$.
Set $V=V(\nu)$ and $\tilde V=V(s_i\nu)$

Let $Z=\{B\in \ev\,;\,\mathop{\oplus}_{\tau\in\om;\inn(\tau)=i}V_{\out(\tau)}
\to V_i\ \text{ is surjective.}\}$. It is clear that $Z$ is an open 
subvariety of $\ev$ and contains both $\CO_{b,\om}$ and $\CO_{b',\om}$.
The group
$GL(V_i)$ acts freely on $Z$.
Let $\pi:Z\to Z/GL(V_i)$ be the projection.

Let $\tilde{Z}=\{\tilde{B}\in E_{\tilde{V},s_i\om}\,;\,\tilde{V}_i\to
\mathop{\oplus}_{\tau\in s_i\om;\out(\tau)=i}\tilde{V}_{\inn(\tau)}
\text{ is injective}\}$ 
and let $\tilde{\pi}:\tilde{Z}\to\tilde{Z}/
GL(\tilde{V}_i)$ be the natural projection.

We define a map $\Xi:Z/GL(V_i)\to \tilde{Z}/GL(\tilde{V}_i)$ as
follows. We fix isomorphisms $V_j\simeq \tilde V_j$ ($j\not=i$).
For $B\in Z$, we take an isomorphism
$\tilde V_i\simeq\Ker\Big(\oplus_{\inn(\tau)=i}V_{\out(\tau)}\to V_i\Big)$.
Then define $\tilde B=\Xi(B)$ by: for $\tau\in s_i\om$,
$\tilde B_\tau$ is $B_\tau$ if $\out(\tau)\not=i$,
and $\tilde B_\tau$ is the composition
$$\tilde V_i\simeq
\Ker\Big(\mathop{\oplus}\limits_{\tau'\in\om;\,\inn(\tau')=i}V_{\out(\tau')}
\to V_i\Big)\hookrightarrow
\mathop{\oplus}\limits_{\tau'\in s_i\om;\,\inn(\tau')=i}
\tilde{V}_{\out(\tau')}\to V_{\out(\bar\tau)}$$
if $\out(\tau)=i$.

It is easy to see that $\Xi$ is well-defined and an isomorphism. 
There are $\tilde{b}$ and $\tilde{b}'\in B(\infty)$ such that
$$
\tilde\pi(\CO_{\tilde b,s_i\om})=
\Xi(\pi(\CO_{b,\om}))\ \text{ and }\ 
\tilde\pi(\CO_{\tilde b',s_i\om})=
\Xi(\pi(\CO_{b',\om})).
$$
Then the equivalence of
(1) and (2) is reduced to
\beq\label{ref}
&&\tilde b=S_i(b)\ \text{ and }\ \tilde b'=S_i(b').
\endeq
In order to see this, set $\bar V=V(\nu+m\alpha_i)$
and take a generic point $B$ of $\Lambda_b$.
Then
\beqn
\bar V_j\simeq
\begin{cases}
V_j&\text{ if $j\not=i$}\\
\Imm\Big(
\mathop\oplus\limits_{\inn(\tau)=i}V_{\out(\tau)}
@>{(B_\tau)}>>V_i@>{(B_\tau')}>>
\mathop\oplus\limits_{\out(\tau')=i}V_{\inn(\tau)}\Big)
&\text{ if $j=i$}\\
\end{cases}
\endeqn
gives a point $\bar B$ of $X_{\bar V}$.
It is easy to see that
$\bar B$ is a generic point of
$\eit^*{}^\max b$ and also
a generic point of
$\eit{}^\max\tilde b$.
Hence we have
$\eit^*{}^\max b=\eit{}^\max \tilde b$
and (\ref{ref}) follows.
\end{pf}
\subsection{}
Only by using Lemma \ref{crt1} and \ref{crt2}
we can show

\begin{prop}
Conjecture 2 is true for $1\leq n \leq 7$.
\end{prop}
In fact we used a computer to check this.

There is a counterexample
in the $n=8$ case derived by the counterexample
in Theorem \ref{counter}.

\begin{example}
Let
\begin{eqnarray*}
w&=&s_1s_3s_2s_4s_3s_5s_4s_3s_2s_1s_6s_7s_6s_5s_4s_3\quad\text{and}\cr
w'&=&s_1s_3s_4s_3s_5s_4s_3s_7.
\end{eqnarray*}
 Here $\{s_i\}_{i\in I}$ are the standard 
generators of symmetric group. Then we have
$$SS(\sideset{^{\pi}}{_{{\cal O}_{w,\om_0}}}{\Bbb C})=
\overline{T^*_{{\cal O}_{w,\om_0}}E_{V,\om_0}}\cup
\overline{T^*_{{\cal O}_{w',\om_0}}E_{V,\om_0}}.$$
This singularity is also realized by a partial flag manifold
as follows.
Let $X'$ be the set of flags $\{F_j\}$ of ${\Bbb C}^8$ with
$0=F_0\subset F_1\subset F_2\subset F_3\subset F_4={\Bbb C}^8$
and $\dim F_j=2j$ ( $j=1,2,3$).
Set $Z=X'\times X'=\{(F,F')\in X'\times X'\}$.
Let $Z_1$ be the $SL(8)$-orbit of $Z$
given by the following table of $\dim \Gr^F_i\Gr^{F'}_j$:

\smallskip
\centerline{
\vbox{\offinterlineskip
\halign
{\strut\hfil\quad$#$\quad&\vrule height16pt#&&\strut\quad\hfil$#$\quad\cr
\hbox{$j\backslash i$}\kern-7pt&height 10pt&1&2&3&4\cr
\noalign{\hrule}
1&&1&0&1&0\cr
2&&0&1&0&1\cr
3&&1&0&1&0\cr
4&&0&1&0&1\cr}}}
\smallskip

and $Z_2$ is given by

\smallskip
\centerline{
\vbox{\offinterlineskip
\halign
{\strut\hfil\quad$#$\quad&\vrule height16pt#&&\strut\quad\hfil$#$\quad\cr
\hbox{$j\backslash i$}\kern-8pt&height 10pt&1&2&3&4\cr
\noalign{\hrule}
1&&2&0&0&0\cr
2&&0&0&2&0\cr
3&&0&2&0&0\cr
4&&0&0&0&2\cr}}}
\medskip
Then $\overline{Y_w}$ (resp. $\overline{Y_{w'}}$)
is the inverse image of
$\overline{Z_1}$ (resp. $\overline{Z_2}$)
by the canonical morphism
$X\times X\to X'\times X'$.
Hence the characteristic variety of
the intersection cohomology sheaf of $Z_1$
contains the conormal bundle of $Z_2$.
The singularity of $Z_1$ at $Z_2$ is the same as
the one of the counterexample in Theorem \ref{counter}.
\end{example}

\end{document}